\documentclass[leqno,12pt]{article}

\usepackage{latexsym}
 \usepackage{graphicx}
\usepackage{amsmath}
\usepackage{amssymb}
\usepackage{mathtools}
\usepackage{algorithm}
\usepackage{algorithmic}
\usepackage{cite}
\usepackage{color}
\usepackage{subfig}

\usepackage[margin=1.3in, top=1.3in, bottom=1.3in]{geometry}

\newcommand{\nc}{\newcommand}
\nc{\nt}{\newtheorem}
\nt{thm}{Theorem}[section]
\nt{cor}[thm]{Corollary}
\nt{prop}[thm]{Proposition}
\nt{lem}[thm]{Lemma}
\nt{defn}[thm]{Definition}
\nt{rem}[thm]{Remark}
\nt{exa}[thm]{Example}
\nt{ass}[thm]{Assumption}
\nt{alg}[thm]{Algorithm}
\nc{\ip}[2]{\mbox{$\langle #1,#2 \rangle$}}
\nc{\pf}{\noindent{\bf Proof\ \ }}
\nc{\finpf}{\hfill{$\Box$}\linespace}
\nc{\linespace}{\vspace{\baselineskip} \noindent}
\nc{\R}{{\bf R}}
\nc{\Rn}{{\bf R}^n}
\nc{\bx}{\bar{x}}
\nc{\by}{\bar{y}}
\nc{\inT}{\mbox{\rm int}\,}
\nc{\cl}{\mbox{\rm cl}\,}
\nc{\gph}{\mbox{\rm gph}\,}
\nc{\argmin}{\mbox{\rm argmin}\,}

\def\tto{\;{\lower 1pt \hbox{$\rightarrow$}}\kern -12pt
           \hbox{\raise 2.8pt \hbox{$\rightarrow$}}\;}
\newenvironment{myequation}{\setcounter{equation}{\value{thm}}
   \begin{equation}}{\addtocounter{thm}{1}\end{equation}}

\nc{\bmye}{\begin{myequation}}
\nc{\emye}{\end{myequation}}

\DeclarePairedDelimiter{\set}{\{}{\}}

\begin{document}
\title{
A simple Newton method for local nonsmooth optimization
}
\author{
A.S. Lewis
\thanks{ORIE, Cornell University, Ithaca, NY 14853, USA.
\texttt{people.orie.cornell.edu/aslewis} \hspace{2mm} \mbox{}
Research supported in part by National Science Foundation Grant DMS-1613996.}
\and
C.J.S. Wylie
\thanks{ORIE, Cornell University, Ithaca, NY 14853, USA.}
}
\date{\today}
\maketitle

\begin{abstract}
Superlinear convergence has been an elusive goal for black-box nonsmooth optimization.  Even in the convex case, the subgradient method is very slow, and while some cutting plane algorithms, including traditional bundle methods, are popular in practice, local convergence is still sluggish.  Faster variants depend either on problem structure or on analyses that elide sequences of ``null'' steps.  Motivated by a semi-structured approach to optimization and the sequential quadratic programming philosophy, we describe a new bundle Newton method that incorporates second-order objective information with the usual linear approximation oracle.  One representative problem class consists of maxima of several smooth functions, individually inaccessible to the oracle.  Given as additional input just the cardinality of the optimal active set, we prove local quadratic convergence.  A simple implementation shows promise on more general functions, both convex and nonconvex, and suggests first-order analogues.
\end{abstract}
\medskip

sequential quadratic programming, nonconvex
\medskip


\section{Introduction}
Accurate deterministic minimization of a smooth objective function $f \colon \Rn \to \R$ constitutes the foundation of the classical optimization literature \cite{nocedal_wright}.  Typically we assume a black-box oracle, returning, at any point $x \in \Rn$, the objective value $f(x)$, gradient $\nabla f(x)$, and perhaps the Hessian $\nabla^2 f(x)$.  Fast local convergence --- superlinear or quadratic --- relies ultimately on Newton's method.

By contrast, nonsmooth optimization, even for convex objectives, is more challenging, at least in the analogous black-box model \cite{nem-yud}.  Armed with the right tools for the circumstances --- some manner of global structural knowledge --- reasonably accurate nonsmooth convex minimization may be tractable.  The contemporary optimization literature highlights various examples: the classical proximal point method and its various splitting-type extensions, when we can suitably decompose the objective into smooth and prox-friendly ingredients \cite{parikh-boyd,ryu-boyd}; basic techniques for smooth minimization (along the lines surveyed in \cite{Nes04b}) extended to nonsmooth objectives amenable to suitable smoothing techniques \cite{nest_nonsm}; polynomial-time interior-point techniques for semidefinite-representable minimization \cite{nesterov-nemirovskii}.

Outside such structured convex optimization realms, however, we fall back on black-box techniques, using function values and subgradients but blind to any global structure beyond perhaps convexity \cite{Nes04b,bubeck}.  The simplicity and all-purpose appeal of such algorithms is tempered by their local convergence in practice.  At large scale, the classical subgradient algorithm is the notoriously slow method of last resort \cite{nem-yud}.  In the convex case, modern cutting plane methods based on analytic centers \cite{goffin-vial,atkinson-vaidya,nesterov-cutting} or the volumetric barrier \cite{vaidya,lee-sidford-wong}, can be highly effective at moderate scale, but convergence is linear.  

Bundle methods, which originate with \cite{lemarechal-bundle,mifflin-semismooth}, have proved quite successful in large-scale practice, are supported by convergence theory \cite{kiwiel,du-ruszcynski}, and  extend to the nonconvex case (although implementation is delicate --- see \cite{kiwiel-book}).  However, classical bundle methods converge sublinearly, and while more recent black-box versions such as \cite{MC05} are faster, they too are complex to implement, and the convergence analysis considers only ``serious'' steps rather than oracle calls.  Sagastiz\'abal's 2018 ICM lecture \cite{claudia-icm} includes an elegant and comprehensive survey.  Impractically large subproblems disadvantage some otherwise promising analogous developments, like level bundle methods \cite{preamb}. 

The quest for superlinear acceleration described in \cite{claudia-icm} has for two decades also driven important ``$\mathcal{VU}$'' and ``partly smooth'' theory, pioneered by Lemar\'echal, Mifflin, Oustry, and Sagastiz\'abal \cite{LS97,LOS,MC02} and paralleled in \cite{Lewis-active}.  However, still missing in black-box nonsmooth optimization, nonconvex or convex, is the fast local convergence of Newton's method.  Our aim here is just such a black-box Newton-type local optimization method, supported by a superlinear convergence theory on a range of interesting nonsmooth functions, and promising in practice.  

As a step in that direction, we take an unconventional {\em semi-structured} approach.  As motivation, we first consider how to minimize a common type of structured nonsmooth objective function, specifically a pointwise maximum of finitely-many smooth component functions, but using a black-box oracle for the objective that {\em cannot access the component functions individually}.  The resulting algorithm becomes a model-based black-box ``bundle Newton'' method (terminology coming from \cite{luksan-vlcek}) that converges locally quadratically on nonsmooth functions of this max-type, but that also shows promise in practice on more general nonsmooth functions, both convex and nonconvex.  

Like any local method, including Newton's method itself, globalizing the algorithm is an important challenge, and a topic for future work, as is the first-order analogue we describe briefly at the conclusion of this work.  Nonetheless, this Newton-type local approach seems a promising new ingredient for nonsmooth optimization.

\section{An algorithm for nonsmooth minimization}
We seek a local method for finding a minimizer $\bx \in \Rn$ of a continuous nonsmooth  objective function $f$.
We assume, around every point in some set, $\mathcal D \subset \Rn$, that the objective $f$ is twice continuously differentiable.   Given any point $s \in \mathcal D$, an oracle returns the value $f(s)$, gradient $\nabla f(s)$, and Hessian $\nabla^2 f(s)$, allowing us to build the corresponding linear and quadratic approximations, $l_s$ and $q_s$, to $f$.  Beyond that, we call on no further information about the function $f$.

Since the objective $f$ is nonsmooth, we cannot typically hope to find a point $s \in \mathcal D$ where the gradient $\nabla f(s)$ is small.  Instead we seek a finite set (or {\em bundle\/}) $S \subset \mathcal D$ with small diameter
\[
\mbox{diam}\, S ~=~ \max \{|s-s'| :s,s' \in S \}
\]
(where $|\cdot|$ denotes the Euclidean norm),
and small {\em optimality measure}
\[
\Theta(S) ~=~ \min \big|\mbox{conv}\big(\nabla f(S)\big)\big|.
\]
We call the elements of the bundle {\em reference points}.  In terms of the simplex
\[
\Delta_S ~=~ \Big\{ \lambda \in \R^S_+ : \sum_{s \in S} \lambda_s = 1 \Big\},
\]
the optimality measure is
\bmye \label{measure}
\Theta(S) ~=~ 
\min_{\lambda \in \Delta_S} \Big| \sum_{s \in S} \lambda_s \nabla f(s) \Big|.
\emye
If the set of gradients $\nabla f(S)$ is affinely independent, then the minimization problem (\ref{measure}) has a unique optimal solution $\lambda$.  This vector, which plays an important role in the algorithm, we call the {\em Lagrange multiplier estimate}.  

We describe an algorithm that iteratively updates a bundle $S$ one reference point at a time, stopping if it encounters a point outside the set $\mathcal D$.  Thus $S$ has a fixed size (cardinality) throughout:  a judicious choice of this integer parameter, denoted $k$, is a crucial feature of the algorithm that distinguishes it from more standard cutting plane and bundle methods.  Too small a value of $k$ causes the method to fail because the optimality measure $\Theta(S)$ remains bounded away from zero;  too large a value may lead to an ill-defined Lagrange multiplier estimate.

\subsubsection*{The role of convexity}
The core algorithm below is aimed at strongly convex objectives $f$, even though its statement makes sense without convexity and our interest in it is purely local.  Indeed, much of the analysis that follows is independent of convexity, and a straightforward modification of the core method results in our culminating nonconvex algorithm.

Nonetheless, simple examples show that the unmodified method can fail for nonconvex objectives.  We present this algorithm first because it is easier to motivate, our aim in Section \ref{motivation}.  Most importantly, the final step in our key quadratic convergence proof (Theorem \ref{convergence}) depends crucially on a convexity argument.  Without more ado, here is an informal description of the core algorithm.

\begin{alg}[$k$-bundle Newton method for strongly convex $f$]
\label{bundle-newton}
{\rm
\begin{algorithmic}
\STATE
\REQUIRE{initial bundle $S \subset \mathcal D$ of size $k$, tolerances 
$\bar\epsilon,\bar\delta \ge 0$};
\FOR{$\mbox{iteration} = 1,2,3,\ldots$}
\FOR{$s \in S$}
\STATE  $l_s(\cdot) = f(s) + \nabla f(s)^T(\cdot - s)$;
\STATE  $q_s(\cdot) = l_s(\cdot) + \frac{1}{2}(\cdot - s)^T \nabla^2f(s)(\cdot - s)$;
\ENDFOR
\STATE  $\delta = \Theta(S)$;
\STATE
choose $\lambda \in \Delta_S$ with $|\sum_{s \in S} \lambda_s \nabla f(s) | = \delta$;
\IF{$\mbox{diam}\,S < \bar\epsilon$ and $\delta < \bar\delta$}
\RETURN{\em Stopped:  nearly optimal};
\ELSE
\STATE
choose $\hat x \in \argmin \big\{ \sum_{s \in S} \lambda_s q_s(x) : x \in \Rn,~ 
l_s(x) ~\mbox{equal for all}~ s \in S \big\}$
\ENDIF
\IF{$\hat x \not\in \mathcal D$}
\RETURN{\em Stopped:  nonsmooth point.};
\ELSE
\STATE
choose $s \in S$ minimizing $\Theta\big((S \setminus \{s\}) \cup \{\hat x\}\big)$;
\STATE  $S = (S \setminus \{s\}) \cup \{\hat x\}$;
\ENDIF
\ENDFOR
\end{algorithmic}
}
\end{alg}

\noindent
Notice that the case $k=1$ is just the classical Newton method.

\section{Motivating the method} \label{motivation}

\subsection{The optimality measure}
Basic aspects of the algorithm coincide with traditional cutting plane methods.  The method aims to construct a sequence of bundles $S$ converging to the minimizer $\bx$ (and hence in particular with diameter converging to zero) in such a way that the optimality measure $\delta = \Theta(S)$
also converges to zero:  in that case we say that the algorithm {\em succeeds}.  In practice we deduce approximate optimality when both these measures are small.

The algorithm computes an optimal vector $\lambda$ (the Lagrange multiplier estimate) for the optimality measure expression (\ref{measure}).  If we define the weighted average of the reference points
\bmye \label{weighted}
\bar s = \sum_{s \in S} \lambda_s s,
\emye 
then $f(\bar s)$ is certainly an upper bound on the optimal value $\min f$.  To obtain an approximate lower bound, we use convexity to note that each linear approximation
\[
l_s(x) = f(s) + \nabla f(s)^T(x-s) \qquad (x \in \Rn),
\]
minorizes $f$, and hence so does their weighted average $\sum_s \lambda_s l_s$.  Denoting by $L$ a Lipschitz constant for $f$ on some convex set containing all the reference points, we have 
$|\nabla f(s)| \le L$ for all $s \in S$.  Then for all points $x \in \Rn$ we have
\begin{eqnarray*}
\lefteqn{
f(\bar s) - L \mbox{diam}(S) ~\le~ f(\bar s) - \sum_{s \in S} \lambda_s L|s-\bar s| 
~\le~ f(\bar s) + \sum_{s \in S} \lambda_s \nabla f(s)^T(\bar s - s)
} \\
& &  ~=~
\sum_{s \in S} \lambda_s l_s(\bar s) ~\le~ \sum_{s \in S} \lambda_s l_s(x) + \delta|x-\bar s|
~\le~ f(x) + \delta|x-\bar s|,
\end{eqnarray*}
so in conclusion we have 
\bmye \label{gap}
\min f ~\le~ f(\bar s) ~\le~ \min\{f+\delta |\cdot - \bar s|\} +  L \mbox{diam}(S).
\emye
Thus if both the diameter of the bundle $S$ and the optimality measure $\delta$ are small, then the objective value at the point $\bar s$ lies in the small interval between the minimum values of the objective function and a slightly perturbed function.  In this sense, the current bundle constitutes an approximate certificate of optimality.  

As a consequence of this argument, success of the algorithm requires a certain lower bound on the bundle size $k$, as we discuss next.

\subsection{A lower bound on bundle size:  Carath\'eodory number}
Given any set $\Gamma \subset \Rn$ containing zero in its convex hull, define the {\em Carath\'eodory number} $\mbox{car}\,\Gamma$, to be the minimum size of a subset whose convex hull contains zero.  By Carath\'eodory's theorem, we see
\[
1 ~\le~ \mbox{\rm car}\,\Gamma ~\le~ 1 + \dim(\mbox{conv}\,\Gamma).
\]
Suppose that the set $\mathcal D$ has full measure.  When $\mathcal D$ is simply the set of all points where $f$ is twice continuously differentiable, this assumption typically holds in practice, and in particular if the objective $f$ is semi-algebraic \cite{Coste-semi}.   Define the limiting gradient set
\bmye \label{limiting-gradient}
\Gamma ~=~ 
\big\{ \lim_r \nabla f(x_r) : \lim_r x_r = \bx,~ x_r \in \mathcal D ~ \mbox{for}~r=1,2,\ldots \big\}.
\emye
Since $f$ is locally Lipschitz, being continuous and convex, we have (see \cite{clarke})
\bmye \label{clarke}
0 \in \partial f(\bx) = \mbox{conv}\, \Gamma.
\emye
For the algorithm to succeed, the optimality measure $\Theta(S)$ must converge to zero for some sequence of bundles $S$ converging to the minimizer $\bx$, so zero is a convex combination of $k$ elements of the set $\Gamma$, and hence the lower bound
\[
k \ge \mbox{\rm car}\,\Gamma
\]
must hold.

\subsection{The active subspace} \label{active}
Having computed the optimality measure $\delta = \Theta(S)$, the method next seeks a new reference point. If the current bundle $S$ is close to the minimizer $\bx$, then the cutting plane model 
$\tilde f \colon \Rn \to \R$ defined by
\[
\tilde f(x) ~=~ \max_{s \in S} l_s(x),
\]
minorizes the objective $f$, and approximates it around $\bx$.  Furthermore, at every point $x$ on the {\em active subspace} where all the linear approximations are equal,
\[
M ~=~ \{x \in \Rn : l_s(x)~ \mbox{all equal for all}~ s \in S\},
\]
the cutting plane model has subdifferential
\[
\partial \tilde f(x) = \mbox{conv}\big( \nabla f(S) \big)
\]
and hence nonsmooth slope (the fastest rate of decrease) equal to $\delta$.  Thus when the optimality measure is small, the cutting plane model is approximately minimized throughout the active subset, so there is where the method seeks  a new reference point.

\subsection{An upper bound on bundle size: affine independence} \label{independence}
For algorithmic stability, the minimizer $\bx$ should be close to the active subspace $M$.  Since the bundle $S$ is close to $\bx$, the values of the linear approximations $l_s(\bx)$ are all close to 
$f(\bx)$.  Equivalently, therefore, we ask that the point 
$\big(\bx,f(\bx)\big) \in \Rn \times \R$ should be close to the affine subspace
\[
\{ (x,t) \in \Rn \times \R : l_s(x) = t ~\mbox{for all}~ s \in S \}.
\]
As we have observed, the residual in the linear system defining this subspace is small at the point $\big(\bx,f(\bx)\big)$. Standard linear algebra shows that this point is therefore close the solution set (which in particular is nonempty) providing that the smallest singular value of the matrix for the system is bounded below by some fixed tolerance 
$\sigma > 0$.  That $(n+1)$-by-$k$ matrix has columns $\binom{\nabla f(s)}{1}$ (for $s \in S$), so this lower bound amounts to uniform affine independence of the gradients. 

Assuming that this condition holds, success of the algorithm then requires an upper bound on the bundle size $k$, since by taking a convergent subsequence of the matrices above we arrive at a limiting matrix with $k$ linearly independent columns of the form $\binom{g}{1}$, where each vector $g$ lies in the limiting gradient set $\Gamma$ in (\ref{limiting-gradient}), and hence in the subdifferential $\partial f(\bx)$.  Thus the function $f$ has at least $k$ affine-independent subgradients at $\bx$, from which we deduce the upper bound
\[
k \le 1+\dim\big(\partial f(\bx)\big).
\]

As we shall see, the uniform affine independence  property holds automatically in our convergence proof for functions of max-type.  For more general versions of the algorithm, however, we verify the property as follows.  For any vectors 
$g_i \in \Rn$ indexed by $i$ in a list $I$ of length $k$, we denote by
\[
\sigma_I\{g_i : i \in I\}
\]
the $k$th largest singular value of a matrix with $k$ columns 
$\binom{g_i}{1}$ (for $i \in I$).  This nonnegative number is zero exactly when the list is affine dependent.  Using this notation, we fix a parameter $\sigma > 0$ at the outset, and add the following check at the beginning of each iteration.
\medskip

\begin{algorithmic}
\IF{$\sigma_S\{\nabla f(s) : s \in S\} < \sigma$}
\RETURN{\em Stopped:  affine dependent gradients};
\ENDIF
\end{algorithmic}

\subsection{Choosing the bundle size} \label{choosing}
To summarize, if the algorithm succeeds, the bundle size $k$ must satisfy both upper and lower bounds involving the Carath\'eodory number of the limiting gradient set and the dimension of the subdifferential:
\[
\mbox{car}\, \Gamma ~\le~ k ~\le~ 1+\dim\big(\partial f(\bx)\big).
\]
In general, these lower and upper bounds may be far apart.  For example, for the Euclidean norm 
$f = |\cdot|$ at the point $\bx = 0$, with the set $\mathcal D = \Rn \setminus \{0\}$, the bounds become
\[
2 \le k \le n+1.
\]
However, in the following case of particular interest to us in this work, the two bounds are equal.

\begin{exa}[Convex max functions] \label{max-function}
{\rm
Consider a nonsmooth function of the form
\[
f(x) ~= \max_{i=1,\ldots,k} f_i(x) \quad (x \in \Rn),
\]
for smooth convex functions $f_i \colon \R^n \to \R$, for $i=1,2,\ldots,k$.  
At the point $\bx \in \Rn$, suppose that the function values $f_i(\bx)$ are all equal, with gradients 
$g_i = \nabla f_i(\bx)$, so we have
\[
\partial f(\bx) ~=~ \mbox{conv}\{ g_i : i=1,2,\ldots,k\}.
\]
Assuming that the list of gradients $\{g_i\}$ is affinely independent, we have 
\[
k ~=~ 1+\dim\big(\partial f(\bx)\big).
\]
Furthermore, assuming that $\bx$ is a minimizer, so $0 \in \partial f(\bx)$, the system
\[
\sum_i \lambda_i g_i = 0, \quad \sum_i \lambda_i = 1, \quad \lambda \in \R^k_+
\]
must then have a unique solution $\hat\lambda \in \R^k$.  Choose ${\mathcal D}$ to be the set of points $x \in \Rn$ for which the maximizing index set $\mbox{argmax}_i\,f_i(x)$ is a singleton, so the limiting gradient set (\ref{limiting-gradient}) is
\[
\Gamma ~=~ \{g_i : i=1,2,\ldots,k\}.
\]
Hence the Carath\'eodory number $\mbox{car}\,\Gamma$ is the number of nonzero components of the vector 
$\hat\lambda$, which is $k$ exactly when $\bx$ is in fact a {\em nondegenerate} minimizer, meaning that zero lies in the relative interior of the subdifferential $\partial f(\bx)$.  
}
\end{exa}

In general, estimating the lower bound on the bundle size, the Carath\'eodory number $\mbox{car}\,\Gamma$, seems challenging in practice.  On the other hand, with respect to the upper bound, global nonsmooth optimization methods --- the various methods we discussed in the introduction, including the subgradient method, the proximal point and proximal gradient methods, other splitting methods, bundle and level bundle methods, along with nonsmooth BFGS \cite{BFGS}, and gradient sampling \cite{rob_grad_samp}, for example ---  typically suggest subdifferential dimension information as they progress.  Given any finite set of points $\Omega \subset \mathcal D$ near the minimizer $\bar x$, we can use equation (\ref{clarke}) to estimate
\[
\partial f (\bar x) ~\approx~ \mbox{conv}\big(\nabla f(\Omega)\big).
\]
The dimension of the set on the right-hand side is the rank of a matrix with columns
$\binom{\nabla f(x)}{1}$ (for $x \in \Omega$).
This suggests that a reasonable estimate of the dimension of $\partial f(\bx)$ is the approximate rank --- the number of singular values larger than some tolerance --- of this same matrix.  This approximate rank then might serve as the integer $k$ in our Newton method.

\subsection{The quadratic subproblem} \label{quadratic-subproblem}
At the end of each iteration we update the bundle $S$ by substituting a new reference point $\hat x \in \Rn$ for that point in $S$ whose deletion minimizes the resulting optimality measure.  The Newtonian flavor of the algorithm arises from the choice of $\hat x \in \Rn$, which solves a simple, linearly-constrained quadratic program that we discuss next.

To succeed, the algorithm must generate bundles $S$ that cluster tightly.  By assumption, the corresponding list of gradients $\nabla f(S)$ is robustly affinely independent, so this gradient information is inconsistent with any smooth model for the objective function $f$:   we instead must seek a simple, well-behaved, nonsmooth model.  We use a max function model, motivated by Example \ref{max-function}, and consider twice continuously differentiable functions $f_s \colon \Rn \to \R$ satisfying 
\[ 
f_s(s) = f(s), \quad \nabla f_s(s) = \nabla f(s), \quad 
\nabla^2 f_s(s) = \nabla^2 f(s), \quad f_s(s') < f(s')
\]
for distinct points $s,s' \in S$.  The precise form of these functions is immaterial to the algorithm.
It is nonetheless reassuring to note that such functions always exist.  We could for example define
\[
f_s(x) ~=~ q_s(x) - \alpha|x-s|^4 \quad (x \in \Rn),
\]
for a sufficiently large constant $\alpha > 0$.  On the other hand, if in fact $f$ is a max-function, we could simply consider each $f_s$ as one of the functions comprising the pointwise maximum.

Now we consider the (unknown) function $\tilde f \colon \Rn \to \R$ defined by
\[
\tilde f(x) = \max_{s \in S} f_s(x)  \quad (x \in \Rn),
\]
as our working model of the function $f$:  it agrees with $f$ up to second order at each of the reference points $s \in S$.  We can minimize this model via the classical nonlinear program
\[
\left\{
\begin{array}{lrcll}
\mbox{minimize}     & t \\
\mbox{subject to}   & f_s(x) - t      & \leq   & 0 & (s \in S) \\
                    & x \in \Rn,         & & t \in \R.
\end{array}
\right.
\]
Since the functions $f_s$ are unknown, we cannot solve this problem exactly.
Instead, we consider a feasible solution $\big(\hat s,f(\hat s)\big)$, for any point $\hat s \in \mbox{conv}\,S$, and follow (loosely) a classical sequential quadratic programming approach to improve it.  We remark that SQP techniques have some history in the nonsmooth optimization literature \cite{miller-malick}.  

A standard SQP approach \cite{nocedal_wright} would proceed in two steps, the first of which estimates the Lagrange multipliers.  Taking, as a first approximation, each inequality constraint to be active, we seek to solve 
\[
\min_{\lambda \in \Delta_S} \Big| \sum_{s \in S} \lambda_s \nabla f_s(\hat s) \Big|
\] 
Approximating the point $\hat s$ by the point $s$ in each summand leads to the optimality measure in the algorithm:
\[
    \Theta(S) = \min_{\lambda \in \Delta_S} \Big| \sum_{s \in S} \lambda_s \nabla f(s) \Big|.
\] 

Fixing the resulting Lagrange multiplier estimate $\lambda$, the Lagrangian for the nonlinear program is
\[
(x,t) \mapsto \sum_{s \in S} \lambda_s f_s(x).  
\]
The second SQP step then aims to reduce its quadratic model at the feasible solution 
$\big(\hat s,f(\hat s)\big)$ over a feasible region defined by linearized constraints.
We approximate the quadratic model by the function
\[
(x,t) \mapsto \sum_{s \in S} \lambda_s q_s(x),  
\]
and we approximate the linearized feasible region, using the active subspace from Section \ref{active}, as
\[
\big\{( x,t) : l_s(x)=t ~\mbox{for all}~s \in S \big\}.
\]
We hence arrive at exactly the quadratic subproblem in the algorithm.  The subproblem is feasible, as we saw in Section \ref{independence}, and bounded below by our assumption of strong convexity.

This loose explanation can be tightened.  In fact, the proof of Theorem \ref{quadratic} will show that, if the point $\tilde x \in \Rn$ minimizes the model $\tilde f$, and the quantity $\nu = \max |S - \tilde x|$ is small, then, under reasonable conditions the solution $\hat x$ of the quadratic subproblem satisfies $|\hat x - \tilde x| = O(\nu^2)$.  In other words, the algorithm computes a good approximation of the minimizer of the model function as the next reference point.

\section{A sequential quadratic programming tool}
For the philosophy underlying the bundle Newton method, the tool we describe in this section is central.  It is a slight variant of a standard sequential quadratic programming technique \cite{robinson-perturbed} --- for completeness, we prove it directly.  Convexity plays no role, throughout this section.

Given functions $h_i \colon \R^n \to \R$, for $i=1,2,\ldots,k$, we consider an equality-constrained optimization problem of the form
\[
(P) \qquad
\left\{
\begin{array}{lrcl}
\mbox{minimize}     & c^T y \\
\mbox{subject to}   & h_i(y)    & =     & 0 \quad (i=1,2,\ldots,k) \\
                    & y         & \in   & \R^n.
\end{array}
\right.
\]
For our purposes, a linear objective function $c^T y$ (for some vector $c \in \Rn$) suffices.  

We say that a point $\bar y \in \R^n$ satisfies the {\em strong second-order sufficient conditions} if $h_i(\bar y) = 0$ for each $i$ (feasibility), each $h_i$ is twice continuously differentiable around 
$\bar y$, the list of constraint gradients
\[
T ~=~ \{ \nabla h_i(\bar y) : i=1,2,\ldots,k \}
\]
is linearly independent, and there exists a Lagrange multiplier vector $\bar \lambda \in \R^k$ (necessarily unique) satisfying
\[
\sum_i \bar\lambda_i \nabla h_i(\bar y) = -c \qquad \mbox{and} \qquad 
\sum_i \bar\lambda_i \nabla^2 h_i(\bar y) ~\mbox{positive definite on} ~ T^{\perp}.
\]
The point $\bar y$ is then a strict local minimizer for the problem $(P)$.

We next consider a Lagrange multiplier estimate $\lambda \in \R^k$ close to $\bar\lambda$.  The traditional approach of sequential quadratic programming linearizes the constraints around a trial point close to the minimizer $\bar y$, and replaces the objective by the corresponding quadratic approximation to the Lagrangian $c^T y + \sum_i \lambda_i h_i(y)$.  Instead, we use a {\em different} reference point 
$y_i \in \R^n$ (near $\bar y$) for each constraint, forming the corresponding linear approximations 
$p_i \colon \R^n \to \R$ defined by
\[
p_i(y) ~=~ h_i(y_i) + \nabla h_i(y_i)^T (y-y_i) \quad (i=1,2,\ldots,k).
\] 
We denote by $Y$ the reference points $(y_1,y_2,\ldots,y_k)$ in the product space $(\R^n)^k$, which is close to 
$\bar Y = (\bar y,\bar y,\ldots,\bar y)$, and consider the following quadratic program, parametrized by $Y$ and $\lambda$:
\[
(QP) \qquad
\left\{
\begin{array}{ll}
\mbox{minimize}     & c^T y + 
\sum_i \lambda_i \big( p_i(y) + \frac{1}{2} (y-y_i)^T \nabla^2 h_i(y_i)(y - y_i) \big) 
\\
\mbox{subject to}   & p_i(y) = 0 \quad (i=1,2,\ldots,k) \\
                    & y \in \R^n.
\end{array}
\right.
\]
The new list of constraint gradients $\{ \nabla h_i(y_i) \}$ is also linearly independent, so any minimizer $y \in \R^n$ for the quadratic program $(QP)$ must satisfy the conditions for a {\em stationary point\/}:
\begin{eqnarray*}
p_i(y) &=& 0 \quad (i=1,2,\ldots,k) \\
c + 
\sum_i \lambda_i \big( \nabla h_i(y_i) +\nabla^2 h_i(y_i)(y - y_i) \big)
& = & \sum_i \mu_i \nabla h_i(y_i)
\end{eqnarray*}
for some multiplier vector $\mu \in \R^k$.  

\begin{thm} \label{sqp}
Consider a point $\bar y \in \R^n$ satisfying the strong second-order sufficient conditions for the problem 
$(P)$.  Then for all $Y \in (\R^n)^k$ near $\bar Y$, and any multiplier vector 
$\lambda = \bar\lambda + O(\|Y-\bar Y\|)$ in $\R^k$, the quadratic program $(QP)$ has a unique stationary point $\hat y = \bar y + O(\|Y-\bar Y\|^2)$, which furthermore is the unique minimizer.
\end{thm}

\pf
We can write the stationary point conditions as a linear system:
\[
\big(M(Y,\lambda)\big)(y,\mu) ~=~ b(Y,\lambda),
\]
for a linear operator $M(Y,\lambda)$ on $\R^n \times \R^k$ and a vector 
$b(Y,\lambda) \in \R^n \times \R^k$, both depending continuously on the parameter $(Y,\lambda)$.
When $(Y,\lambda) = (\bar Y,\bar\lambda)$, the corresponding homogeneous system is 
\begin{eqnarray*}
\nabla h_i(\bar y)^T y &=& 0 \quad (i=1,2,\ldots,k) \\
\sum_i \bar\lambda_i  \nabla^2 h_i(\bar y)y 
& = & \sum_i \mu_i \nabla h_i(\bar y).
\end{eqnarray*}
By the second-order sufficient conditions, this system has only the trivial solution.  Hence the operator $M(\bar Y,\bar\lambda)$ is invertible.

As $\gamma = \|Y-\bar Y\| \to 0$ with $|\lambda - \bar\lambda |  = O(\gamma)$, we have
\[
p_i(\bar y) = O(\gamma^2) \quad (i=1,2,\ldots,k)
\]
and
\begin{eqnarray*}
c + \sum_i \lambda_i \big( \nabla h_i(y_i) +\nabla^2 h_i(y_i)(\bar y - y_i) \big) 
&=&
\sum_i (\lambda_i - \bar \lambda_i) \nabla h_i(\bar y) + O(\gamma^2)  \\
&=&
\sum_i (\lambda_i - \bar \lambda_i) \nabla h_i(y_i) + O(\gamma^2).
\end{eqnarray*}
We deduce 
\[
\big(M(Y,\lambda)\big)(\bar y,\lambda - \bar\lambda) ~-~ b(Y,\lambda)
~=~
O(\gamma^2).
\]
The norm of the inverse of $M(Y,\lambda)$ is uniformly bounded for $(Y,\lambda)$ near 
$(\bar Y,\bar \lambda)$, so
\[
(\bar y,\lambda - \bar\lambda) ~-~ \big(M(Y,\lambda)\big)^{-1}\big( b(Y,\lambda) \big)
~=~
O(\gamma^2).
\]
So there exists a unique stationary point $\hat y = \bar y + O(\|Y-\bar Y\|^2)$.

Providing $\gamma$ is sufficiently small, the point $\hat y$ is in fact a strict minimizer for the quadratic program, since it also satisfies the sufficient condition
that $\sum_i \lambda_i \nabla^2 h_i(y_i)$ is positive definite on the subspace
\[ 
\{
z \in \Rn : 
\nabla h_i(y_i)^T z = 0 ~\mbox{for all}~i
\}.
\]
Otherwise there would exists sequences of points
\[
\{ y_i^r : i=1,2,\ldots,k \}, 
\]
and corresponding sequences of multipliers $\lambda_i^r$
satisfying  $y_i^r \to \bar y$ for each $i$ 
and $\lambda_i^r \to \bar\lambda$ as $r \to \infty$, and unit vectors $z^r \in \Rn$ satisfying
\[
(z^r)^T \Big( \sum_i \lambda_i^r \nabla^2 h_i(y_i^r) \Big) z^r \le 0 \quad \mbox{and} \quad
\nabla h_i(y_i^r)^T z^r = 0 ~\mbox{for all}~i.
\]
In that case, after taking a subsequence, we can suppose that $z^r$ converges to a unit vector $z \in \Rn$ satisfying
\[
z^T \Big( \sum_i \bar\lambda_i \nabla^2 h_i(\bar y) \Big) z \le 0 \quad \mbox{and} \quad
\nabla h_i(\bar y)^T z = 0 ~\mbox{for all}~i,
\]
in contradiction to the second-order sufficient conditions.
\finpf

\section{The max function case} \label{max-function-case}
In this section we analyze carefully how the bundle Newton method, Algorithm \ref{bundle-newton}, behaves when minimizing a max function $f \colon \Rn \to \R$, as in Example \ref{max-function}.  Thus the function $f$ has the form
\bmye \label{representation}
f(x) = \max_{i=1,\ldots,k} f_i(x), \quad (x \in \Rn),
\emye
for some {\em unknown} functions $f_i \colon \R^n \to \R$ that are now assumed to be twice continuously differentiable for $i=1,2,\ldots,k$.  

Like the last section, the development in this section does not rely on convexity, with one key exception.  The proof of the final quadratic convergence result, Theorem \ref{convergence} depends critically on strong convexity.

\bigskip
\noindent
{\bf Note.}
We emphasize a crucial point.  Our interest in max functions is as local models for more general objectives, and as a test bed for the general-purpose Algorithm~\ref{bundle-newton}.  We could easily minimize an explicitly described objective of the form (\ref{representation}) by applying a standard algorithm to the equivalent inequality-constrained optimization problem
\[
{\rm (IP)} \qquad
\left\{
\begin{array}{lrcl}
\mbox{minimize}     & t \\
\mbox{subject to}   & f_i(x) - t    & \le   & 0 \quad (i=1,2,\ldots,k) \\
                    & x \in \Rn     &       & t \in \R.
\end{array}
\right.
\]
However, we seek to minimize the objective function $f$ using an oracle {\em with no access to the individual functions} $f_i$.  
\bigskip

We consider Algorithm \ref{bundle-newton} for the objective $f$, on a neighborhood of a point 
$\bx$.  Corresponding to the fixed (but implicit) representation (\ref{representation}) of $f$, we assume  the strong second-order conditions defined below.
\begin{defn} \label{strong}
{\rm
Given a max function representation of the form (\ref{representation}), we say that a point 
$\bx \in \Rn$ satisfies the {\em strong second-order conditions} when the following properties hold.
\begin{itemize}
\item
{\em Full activity\/}:  the values $f_i(\bar x)$ are equal for all $i$.
\item
{\em Independence\/}:  the gradients $\{ \nabla f_i(\bar x) : i=1,2,\ldots,k \}$
are affinely independent.
\item
{\em Stationarity\/}:  
There exists a Lagrange multiplier vector $\lambda \in \R^k_+$ (necessarily unique) satisfying 
$\sum_i \lambda_i = 1$ and 
$\sum_i \lambda_i \nabla f_i(\bar x) = 0$.
\item
{\em Strict complementarity\/}:
$\lambda_i > 0$ for all $i$.
\item
{\em  Second-order sufficiency\/}:
The Lagrangian $\sum_i \lambda_i \nabla^2 f_i(\bar x)$ is positive definite on the subspace
$
\{
z \in \Rn : 
\nabla f_i(\bar x)^T z ~\mbox{equal for all}~i\}
$.
\end{itemize}
}
\end{defn}

These assumptions are closely related both to problem (IP) and to Example \ref{max-function}.  If we consider the feasible point $\big(\bx,f(\bx)\big)$ for (IP), full activity amounts to all the constraints being active, independence amounts to the usual linear independence constraint qualification, and the remaining three conditions correspond exactly to the analogous conditions for (IP).  Classically, these conditions are sufficient for the point $\big(\bx,f(\bx)\big)$ to be a strict local minimizer for (IP), and hence for $\bx$ to be a strict local minimizer for the max function $f$.  We imposed the first four conditions in Example \ref{max-function}:  there, we referred to the strict complementarity assumption as nondegeneracy, since, assuming the first three conditions, it amounts to 
$0 \in \mbox{ri}\big(\partial f(\bx)\big)$ (even in the nonconvex case).

We refer to the disjoint open sets
\[
\mathcal D_i ~=~ \big\{ x \in \Rn : f_i(x) > f_j(x)~ (j \ne i) \big\},
\]
as {\em activity regions\/}:  the values of the functions $f$ and $f_i$ coincide on $\mathcal D_i$, as do their gradients $\nabla f$ and $\nabla f_i$, and their Hessians $\nabla^2 f$ and $\nabla^2 f_i$.  At any point in the open set
\[
\mathcal D = \bigcup_{i=1}^k \mathcal D_i,
\]
we suppose that an oracle returns the value of $f$ along with its gradient and Hessian.  Our challenge in minimizing $f$ is that we only have access to precise information about the value of each component function $f_i$, its gradient and Hessian, on the region $\mathcal D_i$:  elsewhere in $\Rn$ we only know the bound $f_i \le f$.

The algorithm stops if it encounters a point outside $\mathcal D$. This is a reasonable assumption in practice, since the complement $\mathcal D^c$ is typically a small set.  In particular, around the local minimizer $\bx$, since the gradients $\nabla f_i(\bx)$ are all distinct (being affinely independent), $\mathcal D^c$ is contained in the union of the manifolds 
$(f_i - f_j)^{-1}(0)$ for $i \ne j$.  Thus $\mathcal D$ is a dense open set around $\bx$.

We consider a closed ball $B_{\gamma}(\bx)$ of small radius $\gamma > 0$ around $\bx$.  At the outset of the algorithm, we consider a {\em full} bundle $S \subset B_{\gamma}(\bx)$,
meaning that it contains exactly one reference point in each of the $k$ activity regions $\mathcal D_i$. 
We can therefore write
\bmye \label{set}
S ~=~ \{ x_1,x_2,\ldots,x_k \}, \quad \mbox{where} \quad x_i \in \mathcal D_i \quad i=1,2,\ldots,k.
\emye
We will prove that this property is maintained as the algorithm proceeds.  

We denote by $\bar\sigma$ a certain $k$th largest singular value: 
\bmye \label{singular-value}
\bar\sigma ~=~ \sigma_k
\left( \left[
\begin{array}{cccc}
\nabla f_1(\bx) & \nabla f_2(\bx)   & \cdots    & \nabla f_k(\bx)   \\ 
1               & 1                 & \cdots    & 1
\end{array}
\right] \right).
\emye
The affine independence assumption implies $\bar\sigma > 0$.

\begin{prop} \label{dependence-check}
For any tolerance $\sigma \in (0,\bar\sigma)$, full bundles $S$ close to $\bx$ always satisfy the robust affine independence condition $\sigma_S\big(\nabla f(S)\big) > \sigma$.
\end{prop}

\pf
If the result fails, then for each index $i=1,2,\ldots,k$ there exists a sequence of points $(x_i^r)$ in 
the activity region $\mathcal D_i$ converging to the minimizer $\bx$, such that 
\[
\sigma_k \left( \left[
\begin{array}{cccc}
\nabla f_1(x_1^r) & \nabla f_2(x_2^r)   & \cdots    & \nabla f_k(x_k^r)   \\ 
1               & 1                 & \cdots    & 1
\end{array}
\right] \right)
\leq \sigma < \bar \sigma
\]
for infinitely many $r$.
But that contradicts the continuity of the $k$th largest singular value.
\finpf

We now consider Algorithm \ref{bundle-newton} with the choice of tolerances $\bar\epsilon = 0$,
$\bar\delta = 0$, so that the optimality checks never cause the method to stop.  (The affine independence check discussed in Section \ref{active} will not stop the algorithm either, providing we fix the tolerance $\sigma \in (0,\bar\sigma)$.)  We begin each iteration by forming the corresponding linear and quadratic approximations:
\begin{eqnarray*}
l_i(\cdot) &=& f_i(x_i) + \nabla f_i(x_i)^T(\cdot - x_i),  \\
q_i(\cdot) &=& l_i(\cdot) + \frac{1}{2}(\cdot - x_i)^T \nabla^2f_i(x_i)(\cdot - x_i),
\end{eqnarray*}
and estimate the Lagrange multiplier vector.  We have the following result.  

\begin{prop}
For any small radius $\gamma > 0$, there exists a unique minimizer for the problem
\[
\min_{\lambda \in \R^k} \Big\{ \Big| \sum_i \lambda_i \nabla f_i(x_i) \Big| : \sum_i \lambda_i = 1 \Big\},
\]
and it satisfies $\lambda = \bar\lambda + O(\gamma)$.
\end{prop}

\pf
A vector $\lambda \in \R^k$ solves the problem if and only if there exists a number $\alpha \in \R$ such that
\[
\sum_{i=1}^k \lambda_i = 1 \quad \mbox{and} \quad 
\alpha + \nabla f_j(x_j)^T \Big( \sum_{i=1}^k \lambda_i \nabla f_i(x_i) \Big) = 0  \quad (j=1,2,\ldots,k).
\]
This square linear system is defined by an operator depending smoothly on its parameters, the points $x_i$ (for $i=1,2,\ldots,k$).
In the limit, when $x_i = \bar x$ for each $i$, the system becomes
\[
\sum_{i=1}^k \lambda_i = 1 \quad \mbox{and} \quad 
\alpha + \nabla f_j(\bx)^T \Big( \sum_{i=1}^k \lambda_i \nabla f_i(\bar x) \Big) = 0  \quad (j=1,2,\ldots,k).
\]
The corresponding homogeneous system has only the trivial solution, by affine independence of the set 
$\{ \nabla f_i(\bar x) \}$, so the defining operator is invertible.  The unique solution of this limiting system is clearly $(\bar\lambda, 0)$, and the result now follows.
\finpf

As a consequence of this result and the strict complementarity assumption, the Lagrange multiplier estimates are all positive for small radius $\gamma$.
We next turn to the computation of the new reference point $\hat x$.

Being a fully active strict local minimizer of the inequality-constrained problem (IP), the point $\big(\bx,f(\bx)\big)$ is also a local minimizer for the more restrictive problem
\[
(P') \qquad
\left\{
\begin{array}{lrcl}
\mbox{minimize}     & t \\
\mbox{subject to}   & f_i(x) - t    & = & 0 \quad (i=1,2,\ldots,k) \\
                    & x \in \Rn     &   & t \in \R.
\end{array}
\right.
\]
Furthermore, as is easy to verify, it satisfies the strong second-order sufficient conditions with the same Lagrange multiplier vector $\bar\lambda$.  Hence we can apply Theorem \ref{sqp}.  The corresponding quadratic program is
\[
\left\{
\begin{array}{ll}
\mbox{minimize}     & t + 
\sum_i \lambda_i \big( q_i(x) - t \big) 
\\
\mbox{subject to}   & l_i(x) - t = 0 \quad (i=1,2,\ldots,k) \\
                    & x \in \Rn, ~ t \in \R,
\end{array}
\right.
\]
or equivalently, exactly the quadratic program in Algorithm \ref{bundle-newton}:
\[
\left\{
\begin{array}{ll}
\mbox{minimize}     & \sum_i \lambda_i q_i(x)
\\
\mbox{subject to}   & l_i(x)~ \mbox{equal for all}~ i=1,2,\ldots,k \\
                    & x \in \Rn.
\end{array}
\right.
\]
By Theorem \ref{sqp}, this quadratic program has a unique minimizer $\hat x \in \Rn$, which furthermore satisfies $\hat x = \bar x + O(\gamma^2)$.  

The final step in the iteration substitutes the new reference point $\hat x$ for one of the existing reference points.  Since $\hat x$ is also close to the point $\bx$, this substitution produces another full bundle.  To see this, the following idea is the key tool.  

\begin{prop}
Near $\bx$, the optimality measures of full bundles are always strictly less than those of bundles that are not full.
\end{prop}

\pf
If the result fails, then there exists a sequence of full bundles $S_r$, and a sequence of not full bundles $S'_r$, both shrinking to $\bx$, with 
\[
\Theta(S_r) \ge \Theta(S'_r) \qquad (r=1,2,3,\ldots).
\] 
The left-hand side converges to zero, because
\[
0 ~\in~ \mbox{conv}\{\nabla f_i (\bx) : i=1,2,\ldots,k\},
\]
and hence so does the right-hand side.  After taking a subsequence, we can suppose that there is an index $j$ such that $S'_r \cap \mathcal D_j = \emptyset$, and hence
\[
\liminf_r \Omega(S'_r) ~\ge~ \min | \mbox{conv}\{\nabla f_i (\bx) : i \ne j\}|.
\]
But the right-hand side is strictly positive, because
\[
0 ~\not\in~ \mbox{conv}\{\nabla f_i (\bx) : i \ne j\}.
\]
This contradiction completes the proof.
\finpf

\noindent
Since the reference points stay close to the point $\bx$, we next deduce that we maintain full bundles as the algorithm progresses, as follows.
 
\begin{cor}
For any full bundle $S$ near $\bx$, and any new reference point $\hat x \in \mathcal D$ near $\bx$, there is a unique reference point $s \in S$ minimizing the optimality measure of the new bundle 
$S' = (S \setminus \{s\}) \cup \{\hat x\}$, and $S'$ is then also a full bundle.
\end{cor}

\pf
The previous result shows that when $\hat x$ lies in the activity region $\mathcal D_i$, the unique optimal choice of $s$ is the unique reference point in $\mathcal D_i$.  The result  then follows.
\finpf

In summary, we have proved the following result.  

\begin{thm} \label{quadratic}
Given a max function representation of the objective
\[
f(x) = \max_{i=1,\ldots,k} f_i(x) \qquad (x \in \Rn),
\]
suppose the point $\bx$ satisfies the strong second-order conditions in Definition \ref{strong}.  Then there exists a constant $M>0$ such that the $k$-bundle Newton method (Algorithm \ref{bundle-newton}), with the tolerances $\bar\epsilon = 0$, $\bar\delta = 0$, has the following property.  
Any iteration starting with a full bundle $S$ sufficiently close to 
$\bx$ generates a new reference point $\hat x$ satisfying
\[
|\hat x - \bx| ~\le~ M \max_{s \in S} |s-\bx|^2,
\]
and assuming $\hat x \in \mathcal D$, then
generates a new full bundle by substituting $\hat x$ for the unique reference point in $S$ from the same activity region.
\end{thm}

While this is a suggestive result, it does not yet guarantee convergence.  
To ensure that the sequence of bundles shrinks to $\bar x$, we finally call on our assumption that the objective $f$ is strongly convex.  (In the next section we discuss how to modify the algorithm to handle nonconvex objectives.)

We first develop a simple tool.  As usual, for vectors $z \in \R^k$ we define $\|z\|_{\max} = \max_j |z_j|$.

\begin{lem} \label{tool}
Given constants $\epsilon,M > 0$, consider any sequence of vectors in the orthant $\Rn_+$ with the property that, for each successive pair $z,z'$ in the sequence, there exists an index $i$ such that $z'_j = z_j$ for all $j \ne i$, and furthermore
\[
z_i \ge \epsilon \|z\|_{\max}
\qquad \mbox{and} \qquad
z'_i \le M \|z\|_{\max}^2.
\]
Then providing the initial vector is sufficiently small, the sequence converges to zero at a $k$-step quadratic rate.
\end{lem}

\pf
By induction we see that $\|z\|_{\max}$ is nondecreasing as the vector $z$ evolves along the sequence, providing that the initial vector is sufficiently small. Suppose $z = z_{\mbox{\scriptsize old}}$ at the outset of some iteration, and set 
$\theta = \|z_{\mbox{\scriptsize old}}\|_{\max}$.  At this and every subsequent iteration, the updated component $z_i$ is always set to a new value in the interval $(0,M\theta^2]$.  Each updated component therefore cannot be updated again until we have 
\[
\max_{j=1,2,\ldots,k} z_j ~\le~ \frac{M\theta^2}{\epsilon}.
\]
This inequality must therefore hold after at most $k$ iterations, at which point, 
if $z = z_{\mbox{\scriptsize new}}$, we have
\[
\|z_{\mbox{\scriptsize new}}\|_{\max}
 ~\le~ \frac{M}{\epsilon} \|z_{\mbox{\scriptsize old}}\|^2_{\max},
\]
which completes the proof.
\finpf

We can now prove our main result.

\begin{thm}[Fast convergence for strongly convex max functions]  \label{convergence}
\mbox{}\\
Given a max function representation of the objective
\[
f(x) = \max_{i=1,\ldots,k} f_i(x) \qquad (x \in \Rn),
\]
suppose that the point $\bx$ satisfies the strong second-order conditions in Definition \ref{strong}, with each Hessian $\nabla^2 f_i(\bx)$ positive definite.  Then, given the tolerances 
$\bar\epsilon = 0$, $\bar\delta = 0$, 
the $k$-bundle Newton method (Algorithm~\ref{bundle-newton}) starting from any full bundle in a neighborhood of $\bx$, either stops at a point outside the set $\mathcal D$, or generates a sequence of full bundles that converge $k$-step quadratically to $\bx$.
\end{thm}

\pf
Given any small radius $\gamma > 0$, each function $f_i$ is $\rho$-strongly convex on $B_{\gamma}(\bx)$, for some constant $\rho > 0$.  In fact, the proof that follows takes place entirely in the ball $B_{\gamma}(\bx)$, so we lose no generality in assuming that each function $f_i$ is convex.  

According to Theorem \ref{quadratic}, there exists a constant $M > 0$ such that, for small enough 
$\gamma > 0$, starting from any full bundle in the ball $B_\gamma(\bx)$, Algorithm \ref{bundle-newton} produces a sequence of full bundles
\[
S ~=~ \{ x_1,x_2,\ldots,x_k \}, \quad \mbox{where} \quad 
x_i \in \mathcal D_i \cap B_{\gamma}(\bx) \quad i=1,2,\ldots,k,
\]
and at each iteration replaces a reference point, $x_i$ for some index $i$ with a new reference point
$\hat x \in  D_i \cap B_{\gamma}(\bx)$ which furthermore satisfies
\bmye \label{close}
|\hat x - \bx| ~\le~ M \max_{j=1,\ldots,k} |x_j-\bx|^2.
\emye
By the construction, the new reference point $\hat x$ satisfies
\[
l_j(\hat x) = l_i(\hat x) \quad \mbox{for}~ j=1,2,\ldots,k,
\]
and since the functions $f_i$ are twice continuously differentiable, there exists a constant $R > 0$ such that
\[
f_i(x) ~\le~ l_i(x) + \frac{R}{2}|x-x_i|^2 \quad \mbox{for all}~ x \in B_{\gamma}(\bx).
\]
On the other hand, by strong convexity we have
\[
f_j(x) ~\ge~ l_j(x) + \frac{\rho}{2}|x-x_j|^2 \quad \mbox{for all}~ x \in B_{\gamma}(\bx),~ j=1,2,\ldots,k.
\]
Therefore we deduce
\[
l_i(\hat x) + \frac{R}{2}|\hat x-x_i|^2 
~\ge~ 
f_i(\hat x)
~=~ 
f_j(\hat x)
~\ge~
l_j(\hat x) + \frac{\rho}{2}|\hat x-x_j|^2
~=~
l_i(\hat x) + \frac{\rho}{2}|\hat x-x_j|^2,
\]
and setting $\alpha = \sqrt{\rho/R} > 0$ gives
\[
|\hat x-x_i| ~\ge~ \alpha |\hat x-x_j| \quad \mbox{for}~j=1,2,\ldots,k.
\]

Let $\beta = \max_j |x_j - \bx|$, so by inequality (\ref{close}) we have
$|\hat x - \bx| \le M\beta^2$.  The preceding inequality implies, for each $j=1,2,\ldots,k$, the inequality
\begin{eqnarray*}
\alpha |x_j - \bx| 
&\le& 
\alpha |x_j - \hat x| + \alpha |\hat x - \bx| 
~\le~ 
|x_i - \hat x| + M\alpha\beta^2 \\
&\le& 
|x_i - \bx| + |\bx - \hat x| + M\alpha\beta^2
\le |x_i - \bx| +  M(\alpha+1)\beta^2
\end{eqnarray*}
Maximizing over $j$ implies
\[
|x_i - \bx| \ge \alpha\beta - M(\alpha+1)\beta^2 \ge \frac{\alpha}{2}\beta
\]
providing the radius $\gamma$ is sufficiently small.  So we conclude that the old reference point $x_i$ and the new reference point $\hat x$ satisfy the two key inequalities
\[
|x_i - \bx| \ge \frac{\alpha}{2}\max_j |x_j - \bx|
\qquad \mbox{and} \qquad
|\hat x - \bx| \le M\max_j |x_j - \bx|^2.
\]
We now define
\[
z_j = |x_j - \bx| \in (0,\gamma) \qquad  (j=1,2,\ldots,k), 
\] 
and apply Lemma \ref{tool} to complete the proof.
\finpf  

Notice that the final assumption in the strong second-order conditions --- that the Hessian of the Lagrangian is positive definite on the tangent subspace --- is in fact superfluous here, since we are assuming that each Hessian $\nabla^2 f_i(\bx)$ is positive definite.

\section{Minimizing smooth-nonsmooth sums}
With easy modifications, we can use the same Newton bundle idea as in Algorithm~\ref{bundle-newton} to minimize a function of the form $F = f + r$, where the first component $f \colon \Rn \to \R$ is strongly convex but nonsmooth as before, and the second component $r \colon \Rn \to \R$ is nonconvex but smooth (twice continuously differentiable).  

Using the natural choice of optimality measure for a bundle $S \subset \mathcal D$, namely
\bmye \label{new-measure}
\Theta(S) ~=~ \min \big| \mbox{conv} \big( \nabla F(S) \big) \big|.
\emye
the new algorithm proceeds as follows.

\begin{alg}[$k$-bundle Newton method for $F=f+r$]
\label{bundle-newton2}
{\rm
\begin{algorithmic}
\STATE
\REQUIRE{initial bundle $S \subset \mathcal D$ of size $k$, tolerances 
$\bar\epsilon,\bar\delta \ge 0$};
\FOR{$\mbox{iteration} = 1,2,3,\ldots$}
\FOR{$s \in S$}
\STATE  $l_s(\cdot) = f(s) + \nabla f(s)^T(\cdot - s)$;
\STATE  $q_s(\cdot) = F(s) +  \nabla F(s)^T(\cdot - s) + 
\frac{1}{2}(\cdot - s)^T \nabla^2F(s)(\cdot - s)$;
\ENDFOR
\STATE  $\delta = \Theta(S)$;
\STATE
choose $\lambda \in \Delta_S$ with $|\sum_{s \in S} \lambda_s \nabla F(s) | = \delta$;
\IF{$\mbox{diam}\,S < \bar\epsilon$ and $\delta < \bar\delta$}
\RETURN{\em Stopped:  nearly optimal};
\ENDIF
\STATE  
choose $\hat x \in \mbox{argmin}\big\{ \sum_{s \in S} \lambda_s q_s(x) : x \in \Rn,~ 
l_s(x) ~\mbox{equal for all}~ s \in S \big\}$;
\IF{$\hat x \not\in \mathcal D$}
\RETURN{\em Stopped:  nonsmooth point.};
\ELSE
\STATE
choose $s \in S$ minimizing $\Theta\big((S \setminus \{s\}) \cup \{\hat x\}\big)$;
\STATE  $S = (S \setminus \{s\}) \cup \{\hat x\}$;
\ENDIF
\ENDFOR
\end{algorithmic}
}
\end{alg}

We make some immediate comments in comparison with Algorithm \ref{bundle-newton} for minimizing $f$ alone.  Important to notice is that the affine functions $l_s$ (for $s \in S$), whose equality defines the subproblem constraints, are linear approximations for the component $f$, rather than for the objective $F$:  this subtlety is important for the analysis, which relies on strong convexity of $f$.  In contrast, the functions $q_s$ appearing in the subproblem objective must be quadratic approximations for $F$ to ensure quadratic convergence.  

In this new setting we should consider the possibility of an unbounded quadratic subproblem, since the quadratic approximations are no longer necessarily convex.  For the time being, to emphasize the parallel with Algorithm \ref{bundle-newton}, we omit this check, but include it in the final version (still to come).

As in  Algorithm \ref{bundle-newton}, we might also check affine independence of the gradients:
whether we use the gradients of the component $f$ or the objective $F$ is immaterial, since for bundles of small diameter the two differ by roughly a constant due to the continuity of the gradient of the other component $r$.  When $f$ is a max function, the analogue of Proposition~\ref{dependence-check} still holds:  providing we choose a sufficiently small affine dependence parameter, full bundles close to 
$\bx$ always satisfy the check.  Like the unboundedness check, we omit this check for now, but include it in the final version.

The optimality measure (\ref{new-measure}) is the natural extension of the earlier version, and leads to the precise analogue of the approximate optimality property (\ref{gap}), as well as to our strategy for substituting the new reference point at the end of each iteration.

All of our discussion of the choice of bundle size applies in this new case too, simply replacing the old objective $f$ by the new objective $F$, and noting the simple relationship between the subdifferentials
\[
\partial F(x) = \partial f(x) + \nabla r(x) \qquad (x \in \Rn),
\]
where we understand the left-hand side in the Clarke sense \cite{clarke}.  To motivate the quadratic subproblem, we follow exactly the same argument as in Section \ref{quadratic-subproblem}, but modifying the underlying nonlinear program to
\[
\left\{
\begin{array}{lrcll}
\mbox{minimize}     & r(x) + t \\
\mbox{subject to}   & f_s(x) - t      & \le   & 0 & (s \in S) \\
                    & x \in \Rn,      & & t \in \R.
\end{array}
\right.
\]
We note that the constraints are not changed:  they do not involve the function $r$.  Hence, as we noted, neither does the active subspace in the algorithm.

We can then mimic the analysis of the original bundle-Newton method for max functions.  For notational simplicity, we proved the sequential quadratic programming tool, Theorem \ref{sqp}, for a linear objective, but the case we now need, involving a nonlinear objective, is almost identical.  

We thus arrive at the main analysis, Section \ref{max-function-case}.  Given a representation of the nonsmooth component function $f$ in the objective as
\[
f(x) ~=~ \max_{i=1,\ldots,k} f_i(x) \qquad (x \in \Rn),
\]
we assume that the point $\bx \in \Rn$ satisfies the {\em strong second-order conditions} for the problem of minimizing the objective $f+r$.  By this, we  mean that the feasible point $\big(\bx,f(\bx)\big)$ for the appropriate modification of the problem (IP), namely,
\[
\left\{
\begin{array}{lrcl}
\mbox{minimize}     & r(x) + t \\
\mbox{subject to}   & f_i(x) - t    & \le   & 0 \quad (i=1,2,\ldots,k) \\
                    & x \in \Rn     &       & t \in \R,
\end{array}
\right.
\]
has all constraints active, and satisfies the linear independence constraint qualification, stationarity, strict complementarity, and the second-order sufficient conditions. Once again we note that the constraints are unchanged:  they do not involve the function $r$.  With a virtually identical proof , we arrive at the following generalization of Theorem \ref{convergence}.

\begin{thm}[Quadratic convergence for smooth-nonsmooth sums] \label{smooth-nonsmooth-convergence}
\mbox{} \\ 
Given a max function representation of the nonsmooth function
\[
f(x) = \max_{i=1,\ldots,k} f_i(x) \qquad (x \in \Rn),
\]
and a twice continuously differentiable function $r \colon \Rn \to \R$, suppose that the point 
$\bx \in \Rn$ satisfies the strong second-order conditions for minimizing the objective $f+r$, with each Hessian $\nabla^2 f_i(\bx)$ positive definite.  Then, given the tolerances 
$\bar\epsilon = 0$ and $\bar\delta = 0$, 
the $k$-bundle Newton method (Algorithm~\ref{bundle-newton2}) starting from any full bundle in a neighborhood of $\bx$, either stops at a point outside the set $\mathcal D$, or generates a sequence of full bundles that converge $k$-step quadratically to $\bx$.
\end{thm}

\section{Minimizing weakly convex functions}
Consider a {\em weakly convex} function $F \colon \Rn \to \R$, by which we mean that, for some {\em weak convexity parameter} value
$\eta$, the function $f = F + \frac{\eta}{2}|\cdot|^2$ is convex.  By increasing $\eta$ if necessary, $f$ becomes strongly convex.  Assuming that $F$ is also twice continuously differentiable around every point in the set $\mathcal D$, we can then define the smooth function $r = -\frac{\eta}{2}|\cdot|^2$ and arrive at the following version of Algorithm \ref{bundle-newton2} for minimizing $F$, using the optimality measure:
\[
\Theta(S) ~=~ \min \big| \mbox{conv} \big( \nabla F(S) \big) \big|.
\]

\begin{alg}[$k$-bundle Newton minimization for weakly convex $F$]
\label{bundle-newton3}
{\rm
\begin{algorithmic}
\STATE
\REQUIRE{initial bundle $S \subset \mathcal D$ of size $k$, tolerances 
$\bar\epsilon,\bar\delta \ge 0$, $\sigma > 0$}, \\
weak convexity parameter $\eta$;
\FOR{$\mbox{iteration} = 1,2,3,\ldots$}
\FOR{$s \in S$}
\STATE  $l_s(\cdot) = F(s) + \frac{\eta}{2}|s|^2 + (\nabla F(s)+ \eta s)^T(\cdot - s)$;
\STATE  $q_s(\cdot) = F(s) +  \nabla F(s)^T(\cdot - s) + 
\frac{1}{2}(\cdot - s)^T \nabla^2F(s)(\cdot - s)$;
\ENDFOR
\IF{$\sigma_S\{\nabla F(s) : s \in S\} < \sigma$}
\RETURN{\em Stopped:  affine dependent gradients.}
\ENDIF
\STATE  $\delta = \Theta(S)$;
\STATE
choose $\lambda \in \Delta_S$ with $|\sum_{s \in S} \lambda_s \nabla F(s) | = \delta$;
\IF{$\mbox{diam}\,S < \bar\epsilon$ and $\delta < \bar\delta$}
\RETURN{\em Stopped:  nearly optimal};
\ENDIF
\IF{$\min \big\{ \sum_{s \in S} \lambda_s q_s(x) : x \in \Rn,~ 
l_s(x) ~\mbox{equal for all}~ s \in S \big\} = -\infty$}
\RETURN{\em Stopped:  unbounded subproblem.}
\ELSE
\STATE  
choose optimal $\hat x$;
\IF{$\hat x \not\in \mathcal D$}
\RETURN{\em Stopped:  nonsmooth point.};
\ELSE
\STATE
choose $s \in S$ minimizing $\Theta\big((S \setminus \{s\}) \cup \{\hat x\}\big)$;
\STATE  $S = (S \setminus \{s\}) \cup \{\hat x\}$;
\ENDIF
\ENDIF
\ENDFOR
\end{algorithmic}
}
\end{alg}

We note that the algorithm is almost identical to the original convex version, Algorithm \ref{bundle-newton}, the only change (other than termination checks) being to the definition of the linear approximations defining the active subspace.  As the weak convexity parameter $\eta$ grows large, the active subspace converges to the subspace of all points equidistant from all the reference points in the bundle.

Theorem \ref{smooth-nonsmooth-convergence} then specializes to our culminating result.  Recall that the constant $\bar\sigma > 0$ is defined by equation (\ref{singular-value}).

\begin{cor}[Fast convergence for weakly convex max functions] \mbox{} \\ 
Given a max function representation of the objective
\[
F(x) = \max_{i=1,\ldots,k} f_i(x) \qquad (x \in \Rn),
\]
suppose that the point $\bx \in \Rn$ satisfies the strong second-order conditions given in Definition \ref{strong}.  Then, given the tolerances 
$\bar\epsilon = 0$, $\bar\delta = 0$, small $\sigma > 0$, and sufficiently large weak convexity parameter $\eta$,
the $k$-bundle Newton method (Algorithm~\ref{bundle-newton3}) starting from any full bundle in a neighborhood of $\bx$, either stops at a point outside the set $\mathcal D$, or generates a sequence of full bundles that converge $k$-step quadratically to $\bx$.
\end{cor}

\noindent
{\bf Note:}
Any choice of parameters $\sigma \in (0,\bar\sigma)$ and $\eta$ strictly larger than the largest eigenvalues of each negative Hessian $-\nabla^2 f_i(\bx)$ in fact suffices.
\bigskip

\pf
This follows directly from Theorem \ref{smooth-nonsmooth-convergence}, once we observe that, since each function $f_i$ is twice continuously differentiable, with the given choice of weak convexity parameter 
$\eta$ the functions $f_i + \frac{\eta}{2}|\cdot|^2$ are all strongly convex on a neighborhood of $\bx$, and this local property suffices for the proof.
\finpf

\section{Numerical experiments}

We illustrate the local bundle Newton method on several nonsmooth objective functions.  This handful of simple experiments is meant as a proof of concept rather than comprising any algorithmic recommendations.  Nonetheless, the results appear clearly promising enough to invite future research.

\subsection{Practical considerations} 
We implemented none of the stopping criteria, simply terminating the algorithm manually when rounding error prevented any further progress.

\subsubsection*{Choosing an initial bundle}
In each experiment, we ran a standard global nonsmooth optimization method to
generate a finite set of points $\Omega \subset \mathcal D$ near a minimizer $\bar x$, and used the corresponding gradients to estimate the dimension of the subdifferential $\partial f(\bar x)$ and 
hence choose the bundle size $k$, as we discussed in Section \ref{choosing}.
We then ran a heuristic subset selection procedure \cite{golub-vanloan} to choose a set of
$k$ points in $\Omega$ with robustly affinely independent gradients to form the initial bundle.

For convex problems, we implemented the simple ``Bundle Method with Multiple Cuts'' \cite{du-ruszcynski}, which we reproduce below in our notation.
In the nonconvex case, we implemented the nonsmooth BFGS method \cite{BFGS}.

\begin{alg}[Multiple cut bundle method to minimize convex $f$]
\label{bundle-method}
{\rm
\begin{algorithmic}
\STATE
\REQUIRE{initial bundle $S \subset \R^n$, initial center $z \in S$, stopping tolerance
$\bar\epsilon$, proximal parameter $\rho > 0$, sufficient decrease parameter $\beta \in (0,1)$}
\FOR{$\mbox{iteration} = 1,2,3,\ldots$}
\FOR{$s \in S$}
\STATE $g_s \in \partial f(s)$;
\STATE  $l_s(\cdot) = f(s) + g_s^T(\cdot - s)$;
\ENDFOR
\STATE
Choose $\hat x$ minimizing $\max_{s \in S} l_s(\cdot) + \frac{\rho}{2}|\cdot - z|^2$;
\IF{$f(z) - \max_{s \in S} l_s(\hat x) \leq \bar \epsilon$}
\RETURN{\em Stopped:  nearly optimal.}
\ELSE
\IF{$f(\hat x) \leq f(z) - \beta(f(z) - \max_{s \in S} l_s(\hat x))$}
\STATE $z \leftarrow \hat x$ (serious step)
\ELSE
\STATE
$z \leftarrow z$ (null step)
\ENDIF
\STATE  $S \leftarrow S \cup \{\hat x\}$;
\ENDIF
\ENDFOR
\end{algorithmic}
}
\end{alg}

\noindent
For the bundle method, we chose $\Omega$ to be the set of points whose cutting planes were strongly active in the final iteration.  That is,
\[
    \Omega ~=~ \set{ s \in S : \alpha_s > 0},
\]
where $\alpha_s$ is the dual variable associated with cutting plane $l_s(\cdot)$.
For BFGS, we chose the set $\Omega$ to be the final $2n$ iterates.

\subsubsection*{Solving the quadratic subproblems}
The algorithm involves two quadratic programming subproblems.  
The first involves computing the optimality measure $\Theta(S)$ which amounts to projecting
0 onto the convex hull of vectors $\set{\nabla f(s) : s \in S}$.  
We implemented this as a quadratic program, solved in Gurobi.
For the equality-constrained quadratic programs, 
\begin{equation} \label{primal-qp-subproblem}
\left\{
\begin{array}{lrcl}
\mbox{minimize}     & \sum_s \lambda_s q_s(x)       &       &     \\
\mbox{subject to}   & l_s(x) - t                    & =     & 0 \quad (s \in S) \\
                    & t \in \R,                     &       & x \in \Rn.
\end{array}
\right.
\end{equation}
we simply solve the (linear) optimality conditions directly:
\begin{align} \begin{split} \label{primal-qp-kkt}
  \sum_{s \in S} \lambda_s \nabla^2 f(s)(x - s) + \sum_{s \in S} \mu_s \nabla f(s)
  &= 0, \\
  \sum_{s \in S} \mu_s &= 1, \\
   l_s(x) - t &= 0 \quad (s \in S).
\end{split} \end{align}
The $x$ variable of the solution is then our Newton iterate $\hat x$.

\subsection{Illustrative Examples}
\subsubsection*{A Strongly Convex Problem}

Our first experiment is to minimize max functions of the form
\begin{equation} \label{maxquartfunction}
  f(x) = \max_{i = 1, \ldots k}\set*{g_i^T x + \frac{1}{2} x^T H_i x + \frac{c_i}{24} ||x||^4}
\end{equation}
for $1 \leq k \leq n+1$.
We randomly generate positive constants $c_i$, symmetric positive definite matrices $H_i$, and 
affinely independent random vectors $g_i$ satisfying $\sum_i \lambda_i g_i = 0$
for some $\lambda$ randomly sampled in 
$\set{\lambda > 0 : \sum_i \lambda_i = 1}$.
Then 
\[
    0 \in \partial f(0) = \mbox{conv}\set{g_i : 1 \leq i \leq k}    
\]
so $f$ is nonsmooth at the minimizer of $0$.
The structure of $f$ is unknown to the algorithms, whose access is limited to a 
black box that returns function values, gradients, and Hessians.

\begin{figure}[H]
\centering
\includegraphics[width=0.7\textwidth]{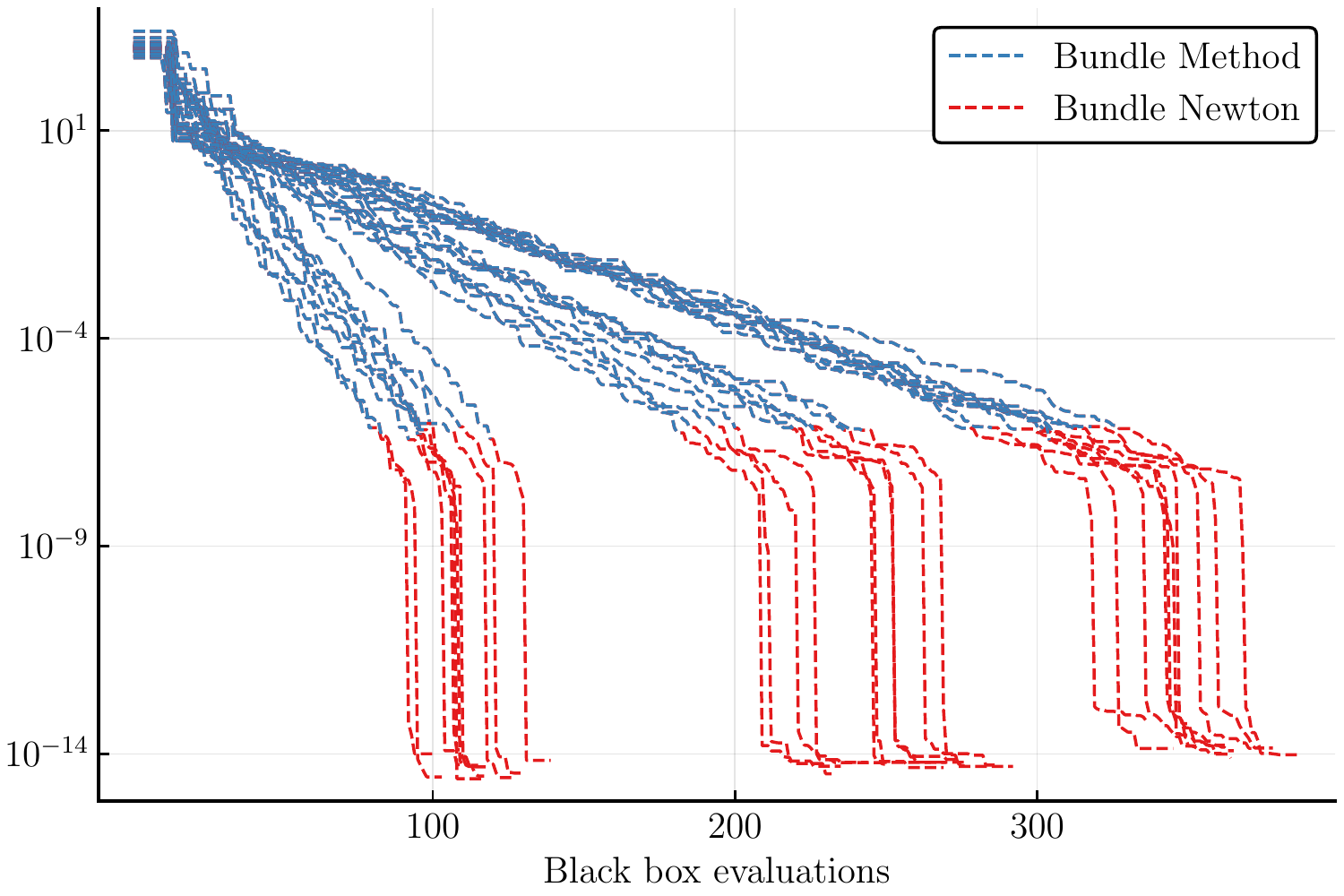}
\caption{Best function value found for the bundle method and
bundle Newton method against number of black box calls for random max functions (\ref{maxquartfunction}) for $k=10,25,40$ in dimension $n=50$.}
\label{maxquart_fevals}
\end{figure}

\begin{figure}[H]
\centering
\subfloat[$\Theta(S)$]{
 \includegraphics[width=0.46\textwidth]{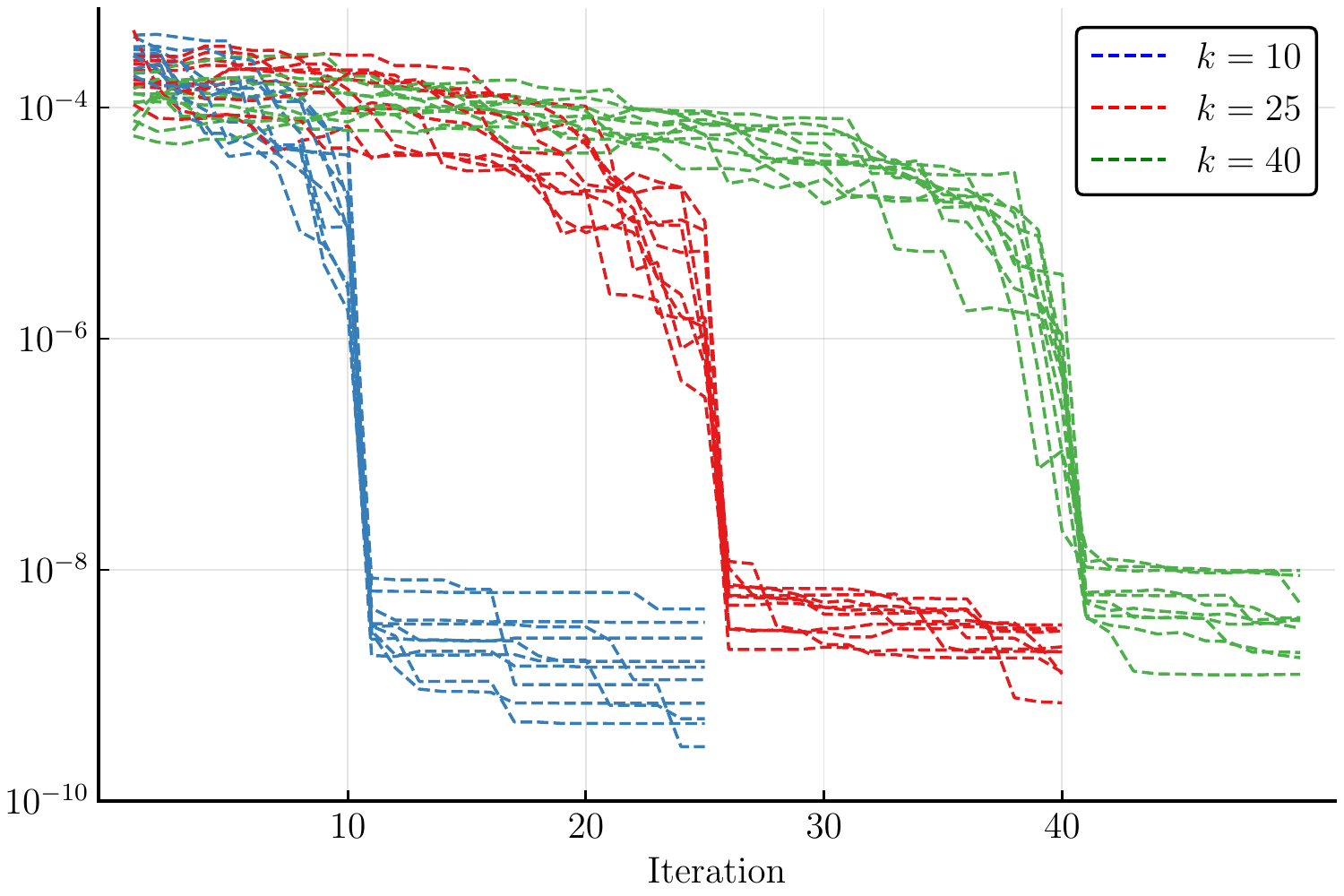}
}
\subfloat[$\mbox{diam}\, S$]{
 \includegraphics[width=0.5\textwidth]{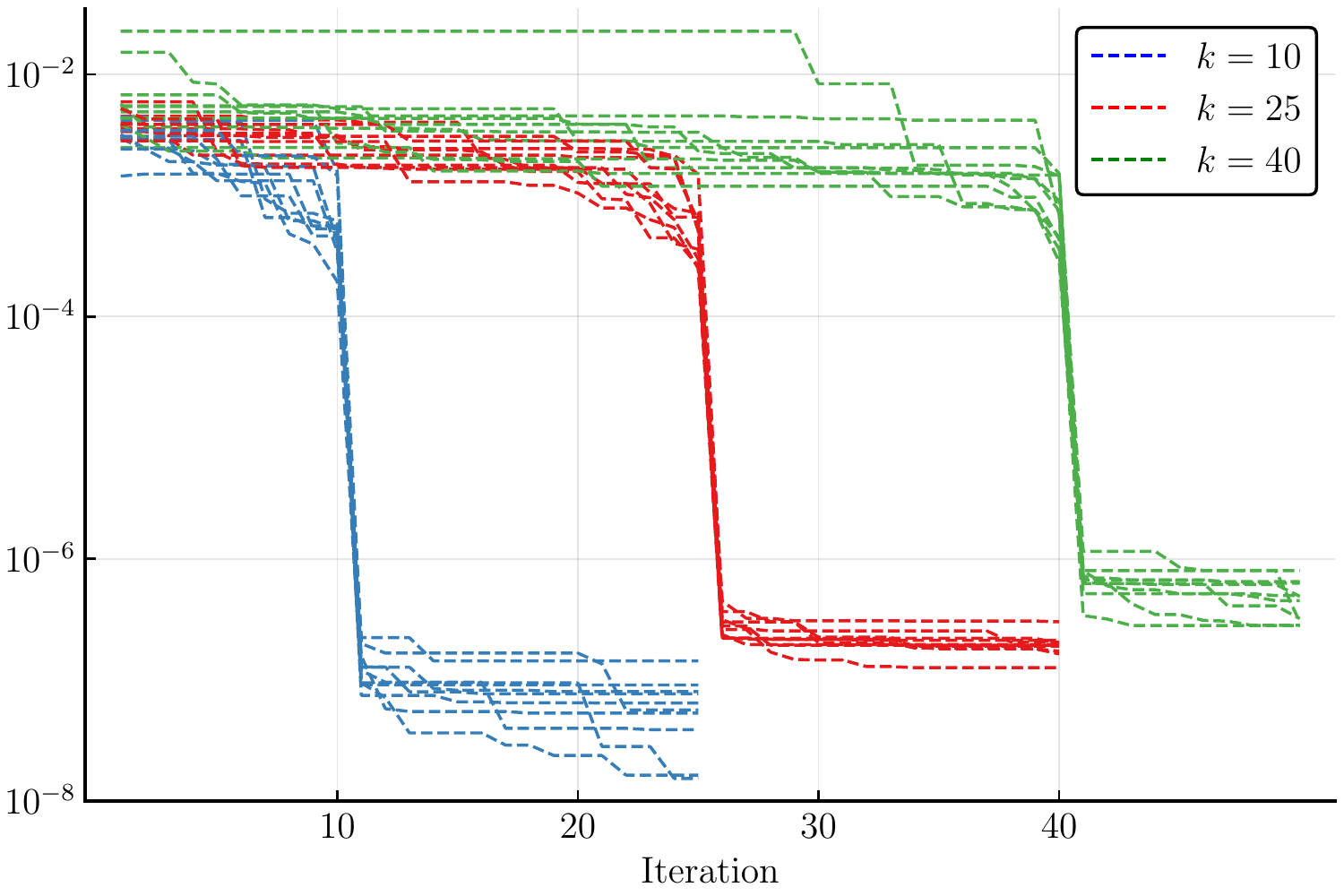}
}
\caption{Optimality measures against iteration count for bundle Newton method for random max functions (\ref{maxquartfunction}).}
\label{maxquart_slope_diameter}
\end{figure}

In random trials for dimension $n=50$, we applied the bundle method, Algorithm~\ref{bundle-method}, in a first phase, with parameters $\rho = 1$ and $\beta = 10^{-5}$ and starting point $z=(1, \ldots, 1)$. The stopping tolerance was set to $10^{-6}$, at which point we initialized and switched to 
the bundle Newton method, Algorithm~\ref{bundle-newton}.  
Results for a number of random trials are shown in Figures~\ref{maxquart_fevals} and \ref{maxquart_slope_diameter}.
For the bundle method phase, we observed a roughly linear rate of convergence of function values to zero, proportional to the degree of nonsmoothness $k$.
Switching to the bundle Newton method results in much more rapid convergence in accordance with the rates predicted by Theorem~\ref{convergence}.

\subsubsection*{A Nonconvex Problem}
To test the nonconvex version of the algorithm, we used Euclidean sum functions of the form
\begin{equation} \label{eucsumfunction}
  f(x) = \sum_{i=1}^k \Big| g_i^T x + \frac{1}{2} x^T H_i x + \frac{c_i}{24} ||x||^4 \Big|
\end{equation}
for $1 \leq k \leq n+1$.  The constants $c_i$, vectors $g_i$ and matrices $H_i$ were randomly generated as in the previous experiment.  As usual, access to $f$ was limited to a black box that returns function values, gradients, and Hessians.

\begin{figure}[H]
\centering
\includegraphics[width=0.63\textwidth]{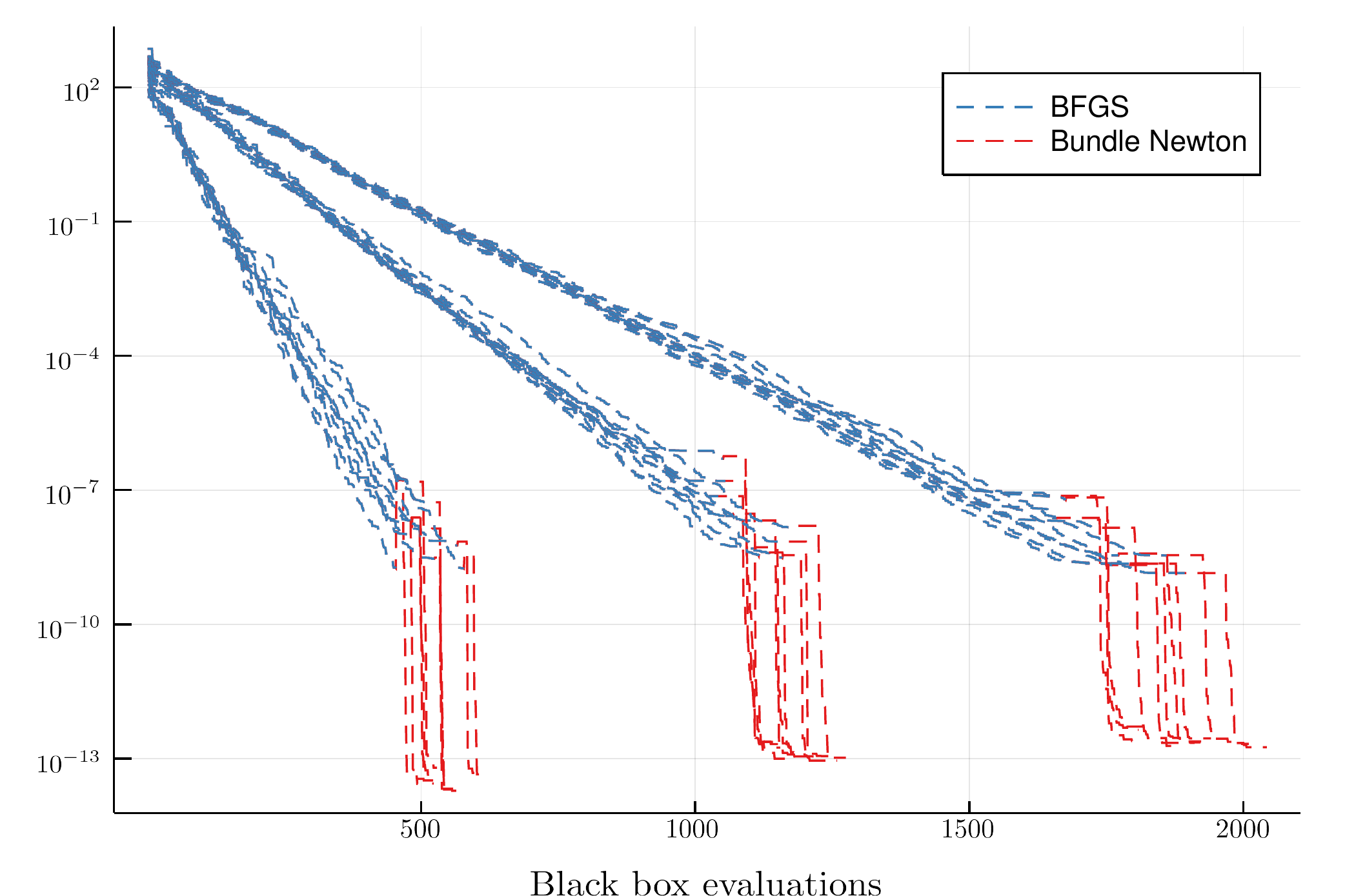}
\caption{Best function value found for BFGS and the
bundle Newton method against number of black box calls on random Euclidean sum functions (\ref{eucsumfunction}) for $k=10, 25, 40$ in dimension $n=50$.}
\label{eucsum_fevals}
\end{figure}

In random trials for dimension $n=50$, we applied nonsmooth BFGS in a first phase
until a breakdown occurred due to numerical instability (as usual with this method \cite{BFGS}).   At this point we switched to the Algorithm~\ref{bundle-newton3} with weak convexity parameter dynamically chosen as
\[
    \eta = \max_{s \in S} \lambda_{\mathrm{max}}\left(-\nabla^2 f(s) \right)
\]
at each iteration.

\begin{figure}[H]
\centering
\subfloat[$\Theta(S)$]{
 \includegraphics[width=0.48\textwidth]{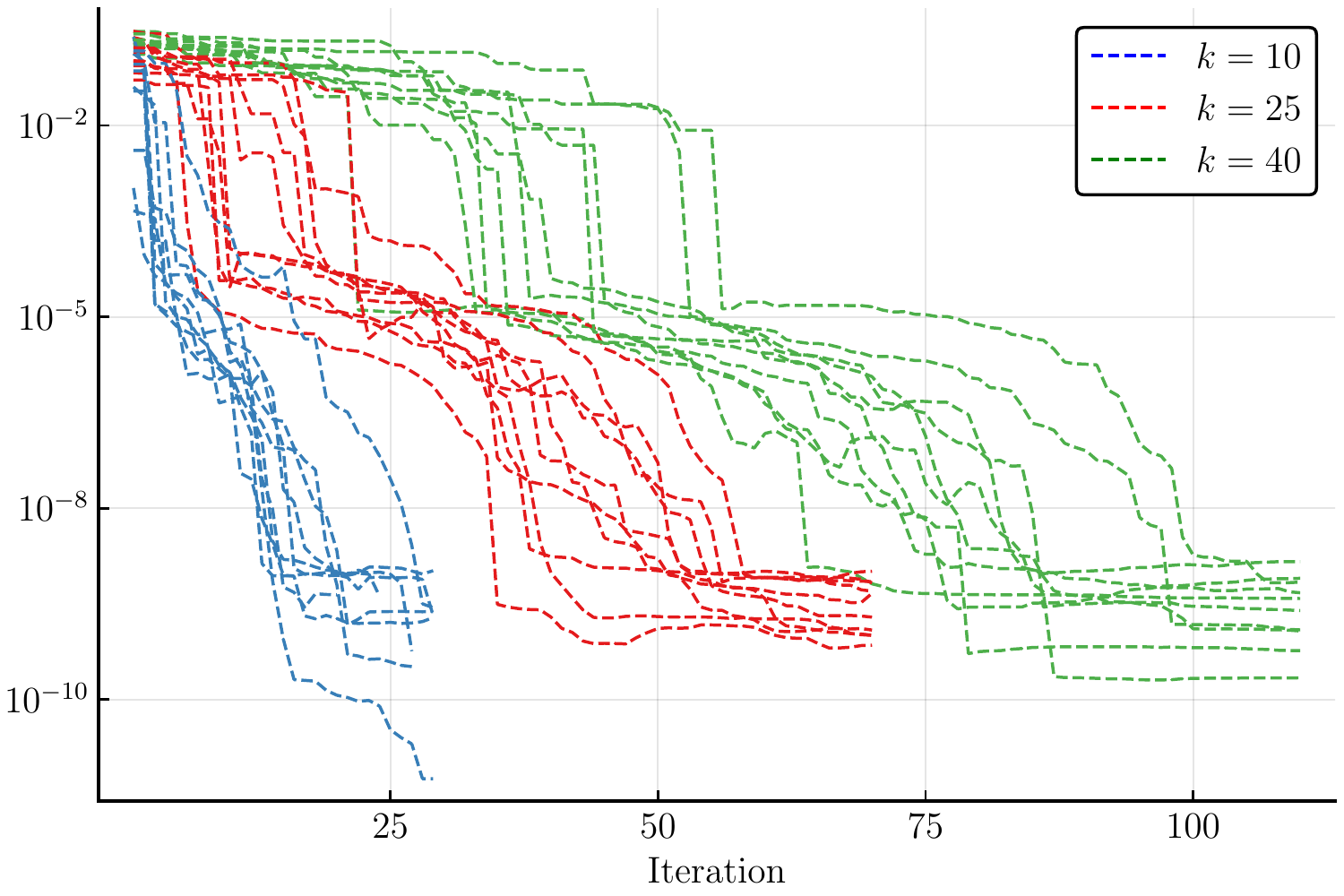}
}
\subfloat[$\mbox{diam}\, S$]{
 \includegraphics[width=0.48\textwidth]{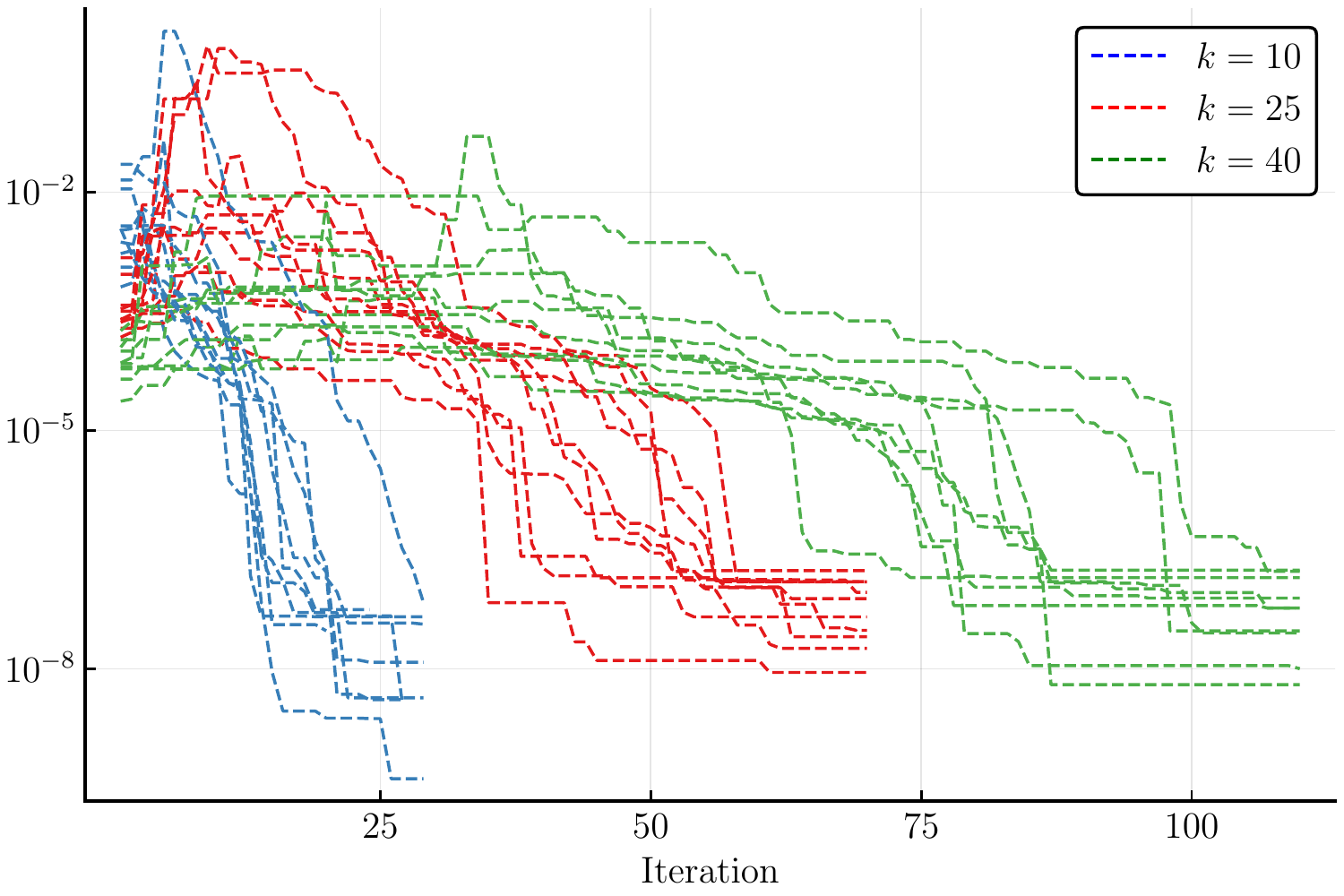}
}
\caption{Optimality measures against iteration count for bundle Newton method for random Euclidean sum functions (\ref{eucsumfunction}).}
\label{eucsum_slope_diameter}
\end{figure}

\subsection{Partly Smooth Functions}

While the theoretical results of this paper were limited to objective functions with finite max structure, experiments suggest that variants of the bundle Newton method may be effective much more broadly.  However, one particular implementation hurdle arises even for simple nonsmooth functions like the Euclidean norm: solving the system 
(\ref{primal-qp-kkt}) directly will be numerically unstable,
due to the ill-conditioning of the Hessians $\nabla^2 f(s)$.

However, for the broad class of \emph{partly smooth} functions in the sense of \cite{Lewis-active}, this ill-conditioning is often highly structured: the nonsmoothness is associated with a certain subspace $\mathcal V$ spanned by the eigenvectors corresponding to large eigenvalues of the Hessian $\nabla^2 f(s)$.  Along an orthogonal manifold, the function behaves smoothly, with well conditioned Hessian.
Motivated by this idea, we follow a simple strategy for solving the system (\ref{primal-qp-kkt}), similar to reduced system approaches for nonlinear programming described in standard texts \cite{nocedal_wright}, and avoiding full Hessian computations.

\subsubsection*{A reduced system approach}

Let $G$ and $b$ be a matrix and vector satisfying
\[
    \set{x : Gx = b} ~=~ \set{x : l_s(x) \text{ equal for all $s \in S$}},
\] 
so that we can write the optimality conditions (\ref{primal-qp-kkt}) as
\begin{align} \begin{split} \label{primal-qp-kkt2}
  \sum_{s \in S} \lambda_s \nabla^2 f(s)(x - s) + G^T \nu
  &= -\sum_{s \in S} \lambda_s \nabla f(s), \\
  Gx &= b.
\end{split} \end{align}
for $\nu \in \R^{k-1}$.
Suppose that we have found matrices $U$ and $V$ such that the matrix
$\begin{bmatrix} U & V \end{bmatrix} \in \R^{n \times n}$ is full rank 
and $GU = 0$ (via a QR factorization of $G^T$, for example).
The columns of $U$ are then a basis for the space $\mbox{Null}(G)$, and 
we can write any solution of (\ref{primal-qp-kkt2}) as $x = U x_u + V x_v$.
The constraint $Gx = b$ then implies $GV x_v = b$,
which can be solved for $x_v$, since we assume $G$ (and hence $GV$) 
is full rank.  We deduce
\[
  \set{x : Gx = b} ~=~ \mbox{Range}(U) + p,
\]
where $p$ is the particular solution $V(GV)^{-1} b$.
Substituting this into the stationarity condition and multiplying through by $U^T$
yields the \emph{reduced} system
\[
  \sum_{s \in S} \lambda_s U^T \nabla^2 f(s)(Ux_u + p  - s) = -\sum_{s \in S} \lambda_s U^T\nabla f(s).
\]
In a slight modification to the algorithm, if we project each reference point onto the active subspace we arrive at the linear system
\[
  \sum_{s \in S} \lambda_s U^T \nabla^2 f(s)U x_u 
  = \sum_{s \in S} \lambda_s \big[ (U^T \nabla^2 f(s)U) U^T(s - p) - U^T\nabla f(s) \big].
\]
This system only involves the \emph{projected Hessians}
$U^T \nabla^2 f(s)U$, which remain well conditioned if the span of $V$ is close to the subspace $\mathcal V$, a property that we have experimentally observed to hold in practice.

\subsubsection*{An Eigenvalue Problem}

Our final experiment to illustrate the reduced systems approach is an eigenvalue problem.  Specifically, given
symmetric matrices $A_0, \ldots, A_n \in \R^{m \times m}$ we seek to minimize
\begin{equation} \label{maxeigfunction}
  f(x) = \lambda_{\max}\Big(A_0 + \sum_{i=1}^n x_i A_i \Big),
\end{equation}
where $\lambda_{\max}(\cdot)$ is the largest eigenvalue function.
Typically minimizers occur at points where $\lambda_{\max}$ has multiplicity $t > 1$, necessitating nonsmooth minimization 
techniques. Under reasonable conditions, the set of points $x \in \Rn$ for which $\lambda_{\max}$ has fixed multiplicity $t$ is a manifold of codimension $\frac{t(t+1)}{2}$, relative to which $f$ is partly smooth \cite{Lewis-active}.

\begin{figure}[H]
\centering
\subfloat[$f(x) - \min f$]{
 \includegraphics[width=0.47\textwidth]{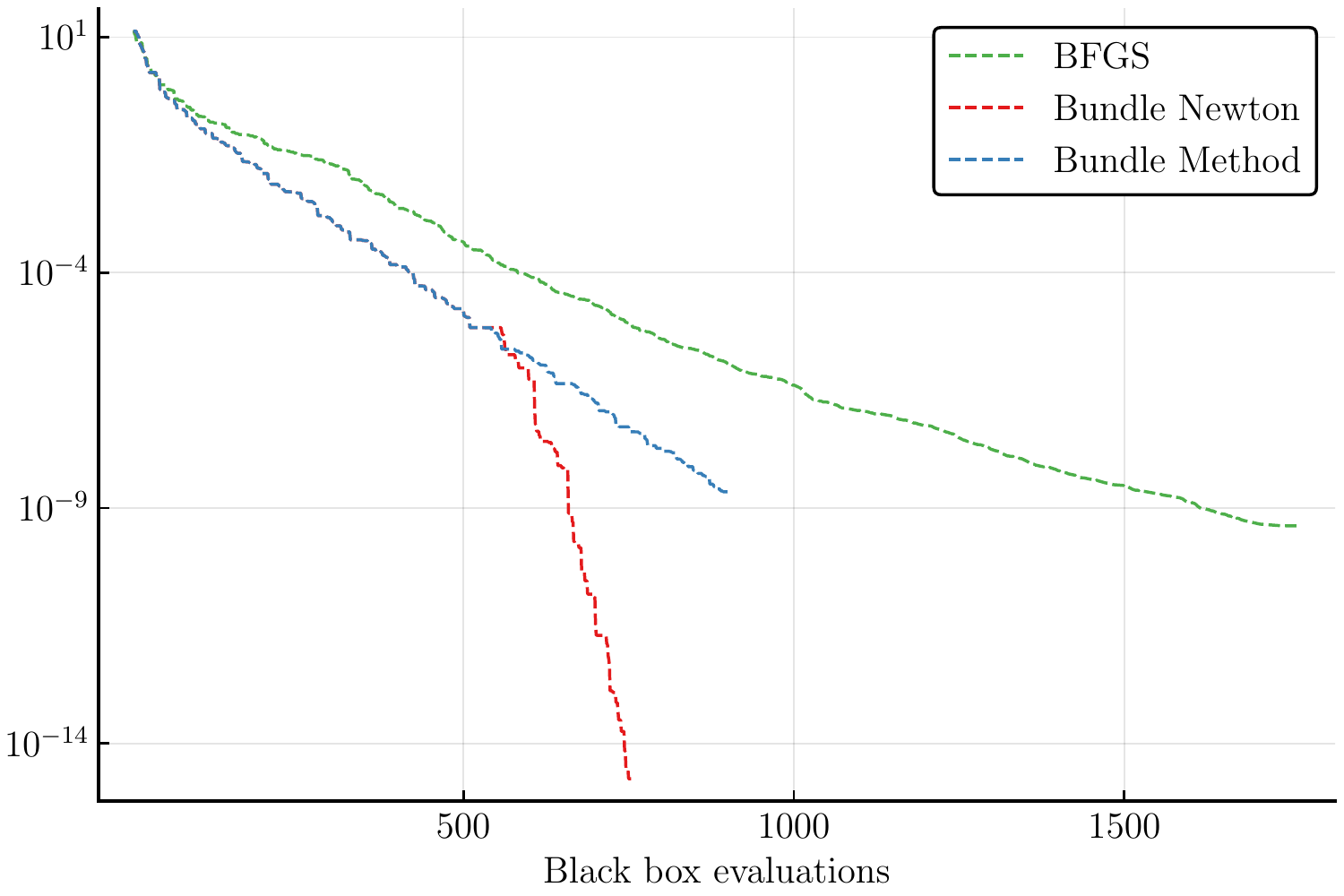}
}
\subfloat[$\mbox{diam}\,S$ and $\Theta(S)$]{
 \includegraphics[width=0.47\textwidth]{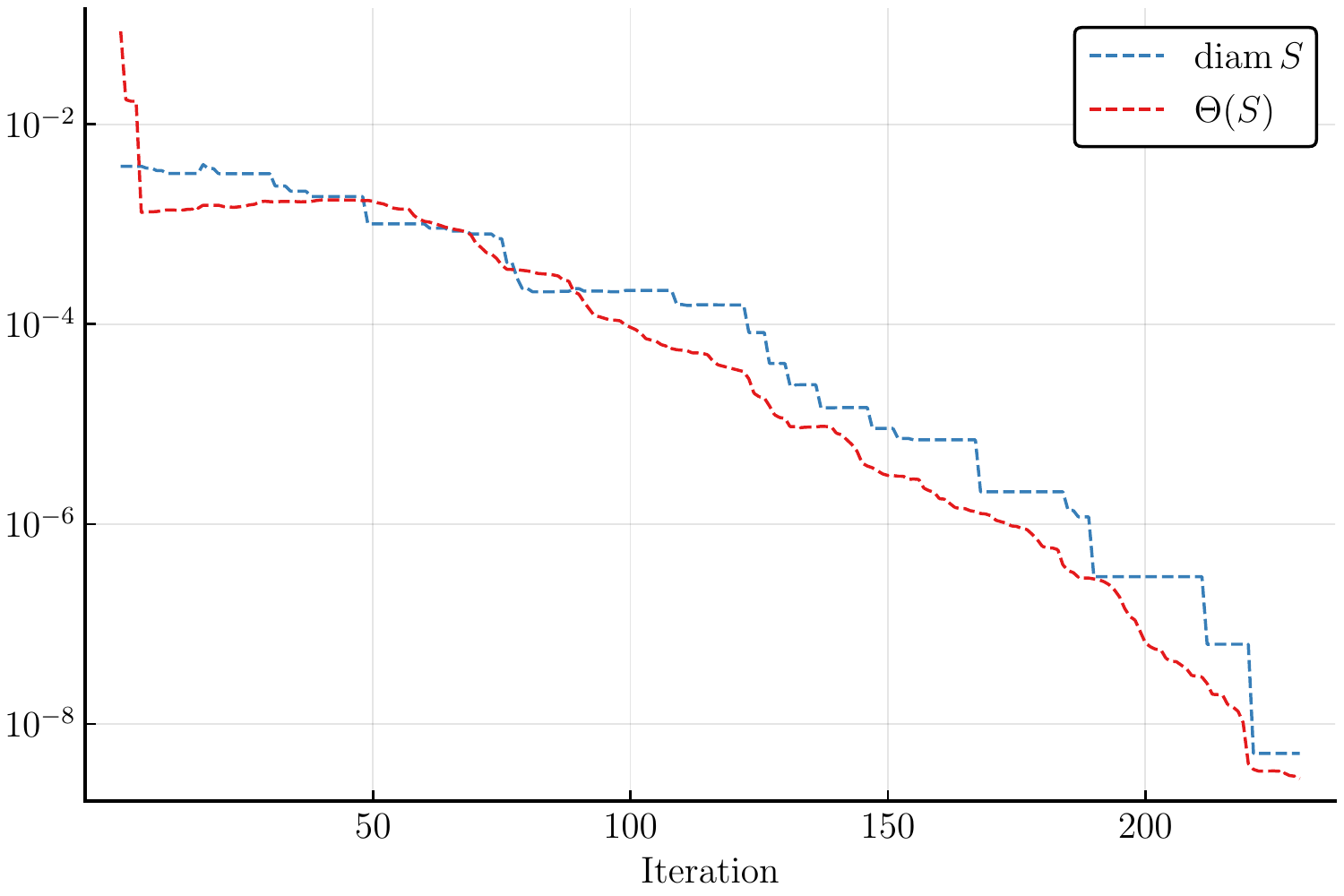}
}
\caption{Function value convergence and optimality measures for the maximum eigenvalue function (\ref{maxeigfunction}) for $n=50$ symmetric matrices in $\R^{25 \times 25}$.}
\label{maxeig_convergence}
\end{figure}

For illustration, in Figure~\ref{maxeig_convergence} we show convergence of the bundle method, BFGS, and bundle Newton method on a typical trial for this problem using random data.  All algorithms were run without termination conditions until numerical issues prevented any further progress.
In this example for $n=50$ matrices in $\R^{25 \times 25}$,
the optimal eigenvalue multiplicity was 6, and we again observe fast 
convergence of the bundle Newton method once the subdifferential dimension 
$\frac{t(t+1)}{2} - 1 = 20$ can be identified.

\begin{figure}[H]
\centering
\subfloat[$\lambda(A(x)) - \min f$]{
 \includegraphics[width=0.47\textwidth]{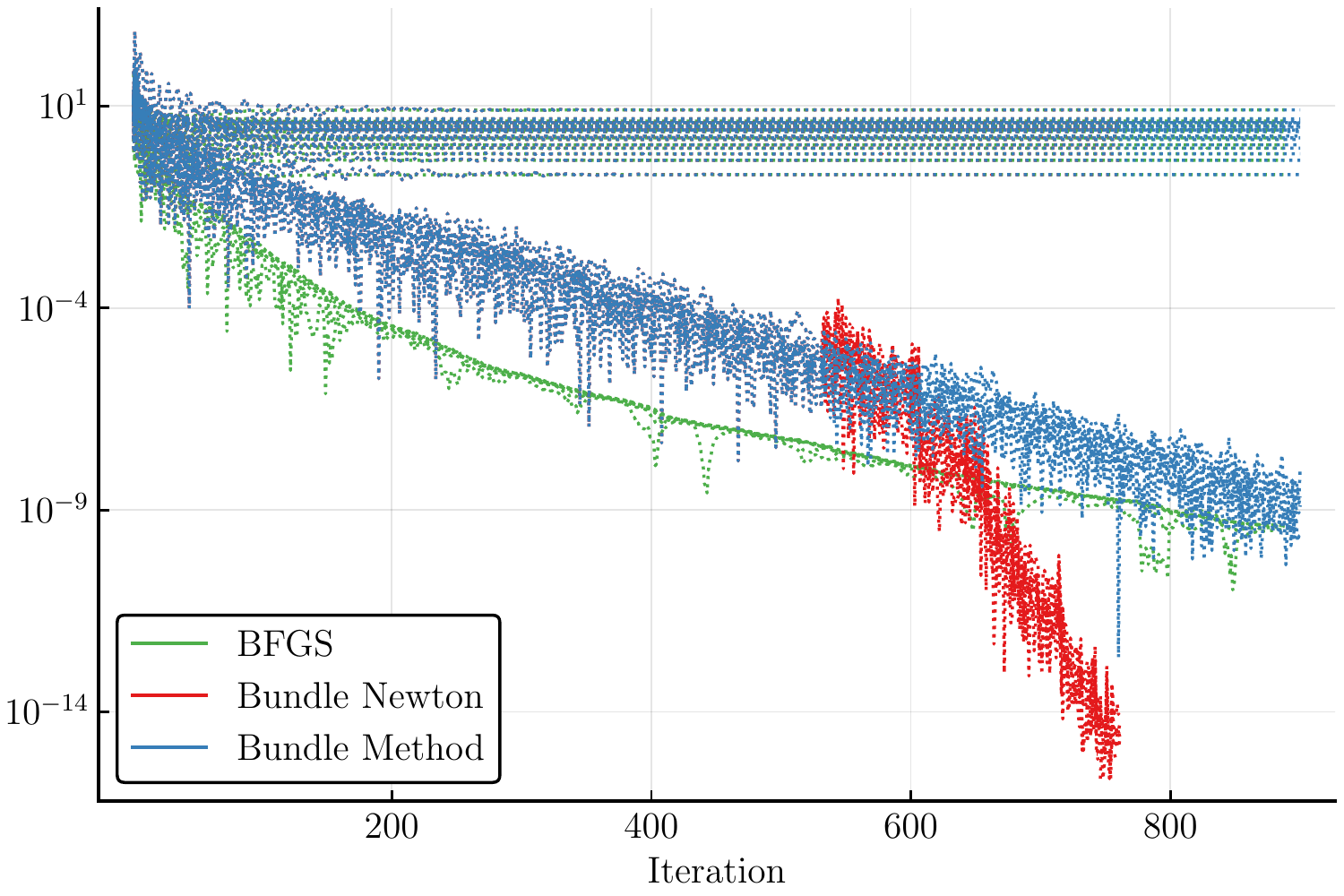}
}
\subfloat[$\lambda_1(A(x)) - \lambda_6(A(x))$]{
 \includegraphics[width=0.47\textwidth]{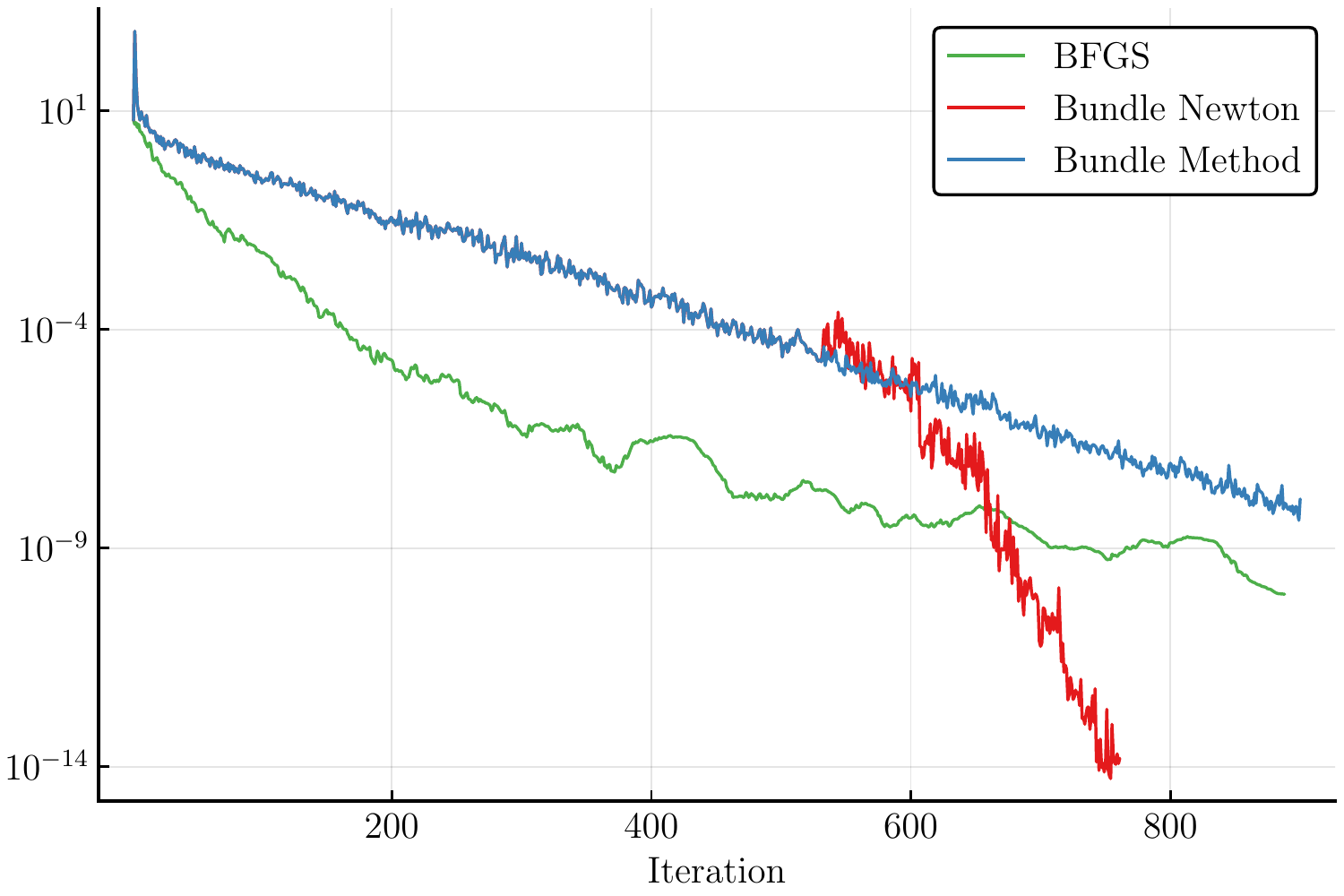}
}
\caption{Clustering of the six largest eigenvalues of $A(x)$.}
\label{maxeig_eigenconvergence}
\end{figure}

(Note that since the optimal objective value is unknown, we instead used the best value found after running the algorithms with a large number of random starting points. This introduces a slight bias in the accuracy reported for the bundle Newton method.)
In Figure~\ref{maxeig_eigenconvergence}, we observe that the bundle Newton method
achieves an eigenvalue clustering several orders of magnitude better than is possible with a bundle method or BFGS.

Using the active manifold to accelerate eigenvalue optimization is not new \cite{overton88,shapiro-fan,oustry99}.  What is remarkable is that the bundle Newton method, combined with a first phase
such as a traditional bundle method, rapidly convergences to the minimizer
\emph{without} any structural knowledge of the function.

\subsubsection*{First-order analogues}
The Newton philosophy that we explore in this work is suggestive even in the more usual setting where Hessians are unavailable.  One straightforward first-order analogue of Algorithm 2,2 replaces the Hessians by suitably tuned multiples of the identity matrix.  Simple implementations seem effective on max functions: a broader investigation is the topic of ongoing work.


\def\cprime{$'$} \def\cprime{$'$}

\end{document}